\documentclass[draft, 11pt]{amsart}
\usepackage{amsxtra}
\usepackage{amssymb}
\usepackage{pb-diagram}
\theoremstyle{plain}
\newtheorem{prop}{Proposition}[section]
\newtheorem{teo}[prop]{Theorem}
\newtheorem{coro}[prop]{Corollary}
\newtheorem{lema}[prop]{Lemma}

\theoremstyle{definition}
\newtheorem{defi}[prop]{Definition}
\newtheorem{ejem}[prop]{Example}
\newtheorem{rem}[prop]{Remark}

\theoremstyle{remark}

\numberwithin{equation}{section}

\newcommand{\N}{\mathbb N}
\newcommand{\Z}{\mathbb Z}

\newcommand{\R}{\mathbb R}
\newcommand{\C}{\mathbb C}

\newcommand{\al}{\alpha}
\newcommand{\g}{\gamma}
\newcommand{\ld}{\lambda}
\newcommand{\Ld}{\Lambda}
\newcommand{\G}{\Gamma}

\newcommand{\f}{\frac}
\newcommand{\tf}{\tfrac}

\newcommand{\arr}{\rightarrow}

\newcommand{\sk}{\smallskip}
\newcommand{\msk}{\medskip}

\newcommand{\vcp}{\Gamma \backslash \R^n}
\newcommand{\I}{\text{\sl Id}}

\newcommand{\man}{M_\Gamma}
\newcommand{\tor}{T_\Ld}
\newcommand{\vep}{\varepsilon}

\newcommand{\on}{\text{O}(n)}
\newcommand{\son}{\text{SO}(n)}
\newcommand{\spin}{\text{Spin}(n)}

\newcommand{\s}{\text{S}}

\newcommand{\hmu}{H_{\mu}^\pm}
\newcommand{\tr}{\text{tr}\,}
\newcommand{\vol}{\text{vol}\,}
\newcommand\dd{\text{d}}

\newcommand{\mjh}{M_{j,h}}

\newcommand{\fa}{\mathcal{F}}

\newcommand{\bjh}{B_{j,h}}

\title[Twisted Dirac operators on compact flat manifolds]{THE SPECTRUM OF TWISTED DIRAC OPERATORS \\ ON COMPACT FLAT MANIFOLDS}
\author[R. J. Miatello]{Roberto J. Miatello}
\address{FaMAF--CIEM \\ Universidad Nacional de C\'ordoba \\ Argentina.}
\email{miatello@mate.uncor.edu, podesta@mate.uncor.edu}
\author[R. A. Podest\'a]{Ricardo A. Podest\'a}

\thanks{Supported by Conicet and Secyt-UNC}

\begin{document}
\bibliographystyle{plain}

\begin{abstract}
Let $M$ be an orientable compact flat Riemannian manifold endowed
with a spin structure. In this paper we determine the spectrum of
Dirac operators acting on smooth sections of twisted spinor
bundles of $M$, and we derive a formula for the corresponding eta
series. In the case of manifolds with holonomy group  $\Z_2^k$, we
give a very simple expression for the multiplicities of
eigenvalues that allows to compute explicitly the $\eta$-series in
terms of values of Riemann-Hurwitz zeta functions, and  the
$\eta$-invariant. We give the dimension of the space of harmonic
spinors and characterize all $\Z_2^k$-manifolds having asymmetric
Dirac spectrum.

Furthermore, we exhibit many examples of Dirac isospectral pairs of $\Z_2^k$-manifolds
which do not satisfy other types of isospectrality.
In one of the main examples, we construct a large family of  Dirac isospectral compact
flat $n$-manifolds, pairwise non-homeomorphic to each other.
\end{abstract}

\maketitle

\section*{Introduction}
 \label{s.intro}
The relation between the geometry of a compact Riemannian manifold $M$ and
the spectral properties of the Laplace operators, $\Delta$ and $\Delta_p$,
acting respectively on smooth functions or on smooth $p$-forms, has been
studied extensively. An elliptic differential operator  whose spectrum is less understood is the Dirac operator $D$. It can  be defined for
Riemannian manifolds having an additional structure, a spin structure. The
goal of the present paper is to investigate properties of the spectrum of Dirac operators on  a compact flat manifold $M$ and with a flat twist bundle. We compare it to other spectral or geometric properties of $M$.

The spectrum of the Dirac operator is explicitly known for a small class
of Riemannian manifolds (see for instance  \cite{Ba} or the list in
\cite{AB}, Table 1). In the context of flat manifolds, Friedrich
(\cite{Fr}) determined the spectrum of $D$ for flat tori
$\tor=\Ld\backslash \R^n$, $\Ld$ a lattice in $\R^n$,  showing the
dependence on the spin structure. In \cite{Pf} Pf\"affle  studied in
detail the Dirac spectrum of 3-dimensional orientable compact flat manifolds,
determining the eta invariants.

The goal of this paper is to study the spectrum of $D$ and the eta
series for an arbitrary compact flat manifold. Such a manifold is
of the form $M_\G:=\G \backslash \R^n$, $\G$ a Bieberbach group.
If $\Ld$ denotes the translation lattice of $\G$, then $F= \Ld
\backslash \G$ is a finite group, the holonomy group of $M_\G$. We
note that by the Cartan-Ambrose-Singer theorem any Riemannian
manifold with finite holonomy group is necessarily flat. As in
\cite{Ch}, we shall use the terminology $F$-manifold for a
Riemannian manifold with finite holonomy group $F$. Already the
case when $F\simeq \Z_2^k$ provides a very large class of
manifolds with a rich combinatorial structure (see Section 3).

The spin structures on $M_\G$ are in one to one correspondence
with homomorphisms $\vep: \G \rightarrow \spin$ satisfying $\mu
\circ \vep =r$ (see Section 1). In \cite{MP} we have shown one can
not hear the existence of spin structures on flat manifolds. For
this purpose, we obtained necessary and sufficient conditions
for the existence of spin structures on $\Z_2^k$-manifolds. 
We note that Vasquez has already shown in \cite{Va} that not every
flat manifold admits a spin structure by giving some examples of
$M$ with holonomy group $\Z_2^k$ having nonzero second
Stiefel-Whitney class, $w_2(M)\ne 0$ (see also \cite{IK} and
\cite{LS}).

Every unitary complex representation  $\rho : \G \rightarrow
\text{U}(V)$, defines a flat vector bundle  $E_\rho:= \G\backslash
(\R^n \times V)$ over $M_\G$. Now, if $(L,\s)$ denotes the spin
representation of $\spin$ we can  consider  the associated twisted
spinor bundle $S_\rho(\man,\vep):= \G\backslash (\R^n\times
(\s\otimes V))$ over $\man$, where $\g\cdot (x, w\otimes v)=(\g x,
L(\vep (\g))(w) \otimes \rho(\g)(v))$, for $w \in \s, v \in V$. We
denote by  $D_\rho$, the  Dirac operator acting on smooth sections
of
$S_\rho(\man,\vep)$. 
For simplicity, we shall assume throughout this paper that $\rho_{|\Ld}=\I$, that is, $\rho$ induces a complex unitary representation of the holonomy group $F\simeq \Ld \backslash \G$. Similar results could be derived in more generality, for instance, assuming that $\rho_{|\Ld}=\chi$, $\chi$ an arbitrary unitary character of $\Ld$. However, this would make the statements of the results more complicated, without gaining too much in exchange.

We now describe the main results. In Theorem \ref{main},   for an
arbitrary compact flat  manifold $\man$, we obtain formulas for
the multiplicities $d_{\rho,\mu}^\pm(\G,\vep)$ of the eigenvalues
$\pm 2\pi \mu$, $\mu > 0$, of the Dirac operator $D_\rho$
associated to a spin structure $\vep$, in terms of the characters
$\chi_{_\rho}$ of $\rho$ and $\chi_{_{L_n}}$,
$\chi_{_{{L_{n-1}}^\pm}}$   of the spin and half spin
representations. The multiplicity formula reads as follows.
If $n$ is odd then
\begin{equation}  \begin{split}
 d_{\rho,\mu}^\pm(\G,\vep)  = \,&
\tfrac{1}{|F|}\,\Big(\!\!\sum_{\scriptsize{\begin{array}{c}\g \in
\Ld \backslash \G\\ B \not\in F_1\end{array}}}
\!\chi_{_\rho}(\g)\; \sum _{u \in (\Ld_{\vep,\mu}^\ast)^B}
e^{-2\pi i u\cdot b}
\;\chi_{_{L_{n-1}^\pm}}(x_\g)\, + \\  \label{multipli}  &    \\
&
\sum_{\scriptsize{\begin{array}{c}\g \in \Ld \backslash \G \\ B
\in F_1\end{array}}} \chi_{_\rho}(\g)\sum _{u \in
(\Ld_{\vep,\mu}^\ast)^B} e^{-2\pi i u\cdot b}
\;\chi_{_{L_{n-1}^{\pm \sigma(u,x_\g)}}}(x_\g)\Big).
\end{split}
\end{equation}
Here $F_1$ is the subset of $F$ corresponding to the elements
$BL_b\in \G$ with $n_B:= \dim\,\text{ker} (B-\I) =1$
and $\Ld^\ast_{\vep,\mu} =\{ u \in \Ld^\ast : \|u\| =\mu, \,  \vep(\ld) =
e^{2\pi i \ld\cdot u} \text{ for every } \ld \in \Ld \}$. Furthermore, for
$\g \in \G$, $x_\g$ is an element in the maximal torus of
$\text{Spin}(n-1)$ conjugate  in $\text{Spin}(n)$ to $\vep(\g)$, and
$\sigma(u, x_\g)$ is a sign depending on $u$ and on the conjugacy class of
$x_\g$ in $\text{Spin}(n-1)$.

If $n$ is even, then the formula reduces to the first summand in (\ref{multipli}) with
$\chi_{_{L_{n-1}^\pm}}$ replaced by $\chi_{_{L_{n-1}}}$.

 We also compute the dimension of the space of harmonic spinors,
(see (\ref{halfharmspinors})), showing that these can only exist
for a special class of spin structures, namely, those restricting
trivially to the lattice of translations $\Ld$.

As a consequence of the theorem, we give an expression, (\ref{etaodd}), for the $\eta$-series $\eta_{(\G,\rho,\vep)}(s)$ corresponding to $D_\rho$ acting on sections of $S_\rho(\man,\vep)$.

In Sections 3 and 4, we restrict ourselves to the case of
$\Z_2^k$-manifolds.
In this case, by computing $\chi_{_{L_{n-1}}}$ and $\chi_{_{L_{n-1}^\pm}}$
we give very explicit expressions 
for the multiplicities. 
Indeed, if the spectrum is  symmetric, we have that $ d_{\rho,\mu}^\pm(\G,\vep) = 2^{m-k-1}\, d_\rho \, |\Ld^\ast_{\vep,\mu}|$ for each $\mu>0$.

The spin $\Z_2^k$-manifolds $(\man,\vep)$ having asymmetric Dirac
 spectrum are of a very special type. This happens 
if and only if $n=4r+3$ and there exists $\g=BL_b\in \G$, with $n_B=1$ and
$\chi_{_\rho}(\g)\not=0$, such that  $B_{|\Ld} =-\delta_\vep \I$.
In this case,  the asymmetric spectrum is the set of eigenvalues
\begin{equation*}
\{\pm 2 \pi \mu_j : \mu_j
=(j+\tfrac12){\|f\|}^{-1},\, j \in \N_0\}
\end{equation*}
where $f$ satisfies $\Ld^B= \Z f$. We furthermore have:
\begin{equation*}
d_{\rho,\mu}^\pm(\G,\vep)  = \left\{\!\! \begin{array}{lc} 2^{m-k-1}\,
\big( d_\rho  \, |\Ld^\ast_{\vep,\mu}| \pm 2 \sigma_\g (-1)^{r+j}
 \chi_{_\rho}(\g) \big) & \mu=\mu_j,
\sk \\ 2^{m-k-1}\, d_\rho \, |\Ld^\ast_{\vep,\mu}| & \mu \ne \mu_j
\end{array} \right.
\end{equation*}
 where $\sigma_\g \in \{\pm1\}$. 
We also give an explicit expression for the eta series:
\begin{equation*}
\eta_{(\G,\rho,\vep)}(s) =  (-1)^{r}\,  \sigma_\g \,\chi_{_\rho}(\g) \,2^{m-k+1}\,\frac{\|   f \|^s}{(4\pi)^s}  \big( \zeta(s, \tfrac 14) - \zeta(s,
\tfrac34)\big)
\end{equation*}
where $\zeta(s,\alpha)=\sum_{j=0}^\infty \tf{1}{(j+\alpha)^s}$
denotes the generalized Riemann-Hurwitz zeta function for
$\alpha\in (0,1]$. From this we obtain that the $\eta$-invariant
of $M_\G$ equals $\eta_\rho=\pm \chi_{_\rho}(\g)\, 2^{[n/2]-k}$,
the sign depending on $\vep$.
This generalizes a result in \cite{Pf} in the case when $n=3$.
We summarize  these results in Theorem \ref{main2} and Proposition \ref{prop.eta2k}.

In Section 4 we compare  Dirac isospectrality with other types of
isospectrality --see Table 1 below--; namely, isospectrality with
respect to the spinor Laplacian
 $\Delta_{s,\rho}:=-D_\rho^2$ and to the Hodge Laplacian on $p$-forms, $\Delta_p$  for
 $0\le p\le n$.
We also look at the length spectrum of $M$ or {\em $[L]$-spectrum}, and
the weak length spectrum of $M$ or {\em $L$-spectrum}, that is, the set of
lengths of closed geodesics counted with and without
multiplicities, respectively.


The information in Examples 4.3, 4.4, 4.5 is collected in Theorem
4.1. We summarize the  results in the following table that shows
the independence of 
Dirac isospectrality from other notions of isospectrality considered.
\renewcommand{\arraystretch}{.5}
\begin{center}
\begin{table}[h]
\caption{Isospectrality}
\begin{tabular}{|c|c|c|c|c|c|c|}
\hline  & & & &  &  & \\
$D_\rho$          & $\Delta_{s,\rho}$ &  $\Delta_p$ $(0\le p\le n)$            & $[L]$     & $L$       & Ex.   & dim  \\ & & & &   & & \\  \hline  & & &   & & &  \\
{\em Yes}    & {\em Yes}           &  {\em No} {\small (generically)} & {\em No}  & {\em No}  & 4.3 (i)     & $n\ge 3$    \\  & &  & & & & \\  \hline  & & &   & & &  \\
{\em Yes}    & {\em Yes}           &  {\em Yes} (if $p$ odd)               & {\em No}  & {\em No}  & 4.3 (iii)   & $n=4t$      \\  & & &  & & & \\  \hline  & &  &  & & &  \\
{\em No}     & {\em Yes}           &  {\em No}              & {\em No}  & {\em No}  & 4.4 (i)    & $n\ge7$       \\  & &  & & & & \\  \hline  & &  &  & & &  \\
{\em Yes/No} & {\em Yes/No}        & {\em Yes}  ($0\le p \le n$)        & {\em Yes} & {\em Yes} & 4.5 (i)     & $n\ge4$       \\  &  & & & & & \\  \hline  & &  &  & & &  \\
{\em Yes/No} & {\em Yes/No}        & {\em Yes}  ($0\le p \le n$)        & {\em No}  & {\em Yes} & 4.5 (ii)    & $n\ge4$       \\  & &  & & & & \\  \hline  
\end{tabular}
\end{table}
\end{center}

Finally, in Example \ref{expfamily}, starting from Hantzsche-Wendt
manifolds (see \cite{MR}), we construct a  large family (of
cardinality depending exponentially on $n$ or $n^2$) of
$\Z_2^{n-1}$\!-manifolds of dimension $2n$, $n$ odd, pairwise
non-homeomorphic,  and Dirac isospectral to each other.

In Section 5, by specializing our formula for $\eta(s)$,  we give
an expression for the eta series and eta invariant for a
$p$-dimensional $\Z_p$-manifold, for each $p=4r+3$  prime. We
obtain an expression for the $\eta$-invariant involving Legendre
symbols and trigonometric sums and give a list of the values  for
$p\le 503$. If $p=3$, our values are in coincidence with those in
\cite{Pf}. Our formulas for the $\eta$-series  are reminiscent of
those obtained in \cite{HZ} to compute the $G$-index of  elliptic
operators for certain low dimensional manifolds. In the case when
$F \simeq \Z_{n}$ in dimension $n=4r+3$ ($n$ not necessarily
prime), an alternative expression for the $\eta$-invariant in
terms of the number of solutions of certain congruences mod$(n)$,
has been given in \cite{SS}, with explicit calculations in the
cases $n=3,7$.

\section{Preliminaries}   \label{s.prelim}

\subsubsection*{Bieberbach manifolds.}
We first review some standard facts on compact flat manifolds (see
\cite{Ch} or \cite{Wo}). A {\it Bieberbach group} is a discrete,
cocompact, torsion-free  subgroup $\G$ of the isometry group $I(\R^n)$ of
$\R^n$. Such $\G$ acts properly discontinuously on $\R^n$, thus $M_\G =
\G\backslash\R^n$ is a compact flat Riemannian manifold with fundamental
group $\G$. Any such manifold arises in this way and will be referred to as
 a {\em Bieberbach manifold}. Any element $\g \in I(\R^n)=\on \rtimes
\R^n$ decomposes uniquely as $\g = B L_b$, with $B \in \on$ and $b\in
\R^n$. The translations in $\G$ form a normal maximal abelian subgroup of
finite index $L_\Ld$ where  $\Ld$ is a lattice in $\R^n$ which is
$B$-stable for each $BL_b \in \G$. The restriction to $\G$ of the
canonical projection $r:\text{I}(\R^n) \rightarrow \on$ given by
$BL_b\mapsto B$ is a homomorphism with kernel $L_{\Ld}$ and $F:=r(\G)$ is a
finite subgroup of $\on$.
 Thus, given a Bieberbach group $\G$, there is an exact sequence
$0 \arr \Ld \arr \G \stackrel{r}{\arr} F \arr 1$. The group $F
\simeq \Ld \backslash \G$ is called the {\em holonomy group} of
$\G$ and gives the linear holonomy group of the Riemannian
manifold $M_\G$. Since we will be working with spin manifolds, we
shall assume throughout this paper that $M_\G$ is orientable,
i.e.\@ $F \subset \son$. The action by conjugation of $\Ld
\backslash \G$ on $\Ld$ defines an integral representation of $F$,
called the {\em holonomy representation}.
This representation can be rather complicated. For instance,
already in the case when $F \simeq \Z_2^2$, there are Bieberbach
groups with indecomposable holonomy representations for
arbitrarily large $n$.
By an $F$-manifold we understand a Riemannian manifold with holonomy group $F$.
In Sections 3 and 4 of this paper we will study in detail the Dirac spectrum of $\Z_2^k$-manifolds.

Let   $\Ld^\ast= \{\ld'\in \R^n : \ld \cdot \ld'\in \Z \text{ for
any } \ld \in \Ld \}$  denote the dual lattice of $\Ld$ and, for
any $\mu \ge 0$, let $\Ld^*_\mu=\{\ld \in \Ld^* : \|\ld\|=\mu\}$.
We note that this $\Ld^*_\mu$ equals  $\Ld^*_{\mu^2}$ in the
notation used in  \cite{MR4}, \cite{MR5}, in the  study of the
spectrum of Laplace operators.

If $B \in \text{O}(n)$ set
\begin{equation}\label{nB}   \begin{array}{c}
(\Ld^*)^B=\{\ld \in  \Ld^* : B\ld =\ld\},  \qquad   (\Ld^*_{\mu})^B =
 \Ld^*_\mu \cap(\Ld^*)^B, \\ \\
 n_B:= \text{dim}\ker(B-\I)=\dim (\R^n)^B. \end{array}
\end{equation}
 If  $\G$ is a Bieberbach group then the
torsion free condition implies that  $n_B>0$ for any  $\g=BL_b \in \G$.
Since $B$ preserves $\Ld$ and $\Ld^\ast$ we also have that $(\Ld^\ast)^B
\ne 0$.
For such $\G$, we set
\begin{equation}\label{F1}
 F_1 = F_1(\G) : = \{ B \in F=r(\G) : n_B =1 \}.
\end{equation}
where $n_B$ is as in (\ref{nB}). 

\subsubsection*{Spin group}
Let $Cl(n)$ denote the Clifford algebra of $\R^n$ with respect to
the standard inner product $\langle \,,\, \rangle$ on $\R^n$ and
let $\C l(n)=Cl(n)\otimes \C$ be its complexification. If $\{e_1,
\dots, e_n\}$ is the canonical basis of $\R^n$ then a basis for
$Cl(n)$ is given by the set $\{ e_{i_1}\ldots e_{i_k} \,:\, 1 \leq
i_1<\dots <i_k \leq n\}$. One has that $vw+wv+2\langle
v,w\rangle=0$ holds for all $v,w\in \R^n$, thus $e_i e_j=-e_je_i$
and ${e_i}^2=-1$ for
$i,j=1,\ldots,n$. 
Inside the group of units of
$Cl(n)$ we have the 
spin group given by
\begin{equation*}
\text{Spin}(n)= \{g=v_1 \ldots v_{k} \,:\, \|v_j\| =1
,\;j=1,\ldots,k\text{ even} \}. 
\end{equation*}
which is a compact, simply connected Lie group if $n\geq3$.
 There is a Lie group epimorphism with kernel $\{\pm 1\}$ given by
\begin{equation}\label{mucovering}
\mu:\text{Spin}(n)\rightarrow \text{SO}(n), \qquad v\mapsto (x \rightarrow
vxv^{-1}).
\end{equation}

If $B_j$ is a matrix for $1\le j \le m$, we will abuse notation by
denoting by $\text{diag}(B_1,\dots, B_m)$ the matrix having  $B_j$ in the
diagonal position $j$.
Let $ B(t) = \left[ \begin{smallmatrix}     \cos t & -\sin t \\
    \sin t & \cos t   \end{smallmatrix}\right] $ with $t \in \R$ and for $t_1,\ldots,t_m \in \R$ let
\begin{equation}\begin{split}\label{torielements}
& x_0(t_1,\ldots,t_m) :=\left\{\begin{array}{ll}diag(B(t_1),\ldots,
B(t_m)),
\quad & \text{if }n=2m \sk \\
diag( B(t_1),\dots, B(t_m),1), \quad & \text{if } n=2m+1
\end{array} \right.         \\
& x(t_1,\dots,t_m) := \prod_{j=1}^m (\cos t_j + \sin t_j e_{2j-1}e_{2j}) \in \text{Spin}(n).
\end{split}
\end{equation}
 Maximal tori in   $\text{SO}(n)$ and $\text{Spin}(n)$ are respectively
given by
\begin{equation}\label{tori}
  {T}_0=\Big\{x_0(t_1,\dots,t_m) \;:\; t_j \in \R \Big\}, \quad
T=\left\{x(t_1,\ldots,t_m) :t_j \in \R \right\}.
\end{equation}
The restriction
$\mu:{T}\rightarrow T_0$ is a 2-fold cover (see \cite{LM}) and
$$
\mu(x(t_1,\dots,t_m))
= x_0(2t_1,\ldots,2t_m).
$$
The Lie algebras of $\text{Spin}(n)$ and that of  ${T}$ are
$\mathfrak{g} = \text{span} \{e_i e_j : 1\le i<j\le n\}$ and
$\mathfrak{t}=\text{span} \{e_{2j-1}e_{2j} : 1\le j \le m\}$,
respectively.

In the Appendix we have collected some specific facts on spin groups and spin representations that are used in the body of the paper.

\subsubsection*{Spin structures on flat manifolds}
If $(M,g)$ is an orientable Riemannian manifold of dimension $n$, let
$\text{B}(M)=\bigcup_{x\in M}\text{B}_{x}(M)$ be the bundle of oriented
frames on $M$ and $\pi:\text{B}(M)\rightarrow M$ the canonical projection,
that is, for $x\in M$, $\text{B}_{x}(M)$ is the set of ordered oriented
orthonormal bases $(v_1,\ldots,v_n)$ of ${T}_x(M)$ and
$\pi((v_1,\ldots,v_n))=x$. $\text{B}(M)$ is a principal
$\text{SO}(n)$-bundle over $M$. A {\it spin structure} on $M$ is an
equivariant 2-fold cover $p:\tilde{\text{B}}(M)\rightarrow \text{B}(M)$
where $\tilde{\pi}:\tilde{\text{B}}(M)\rightarrow M$ is a principal
$\text{Spin}(n)$-bundle and $\pi \circ p=\tilde{\pi}$. A manifold in which
a spin structure has been chosen is called a {\em spin manifold}.

We will be interested on spin structures on quotients
$\man=\Gamma\backslash \R^n$,  $\Gamma$ a Bieberbach group. If $M=\R^n$,
we have $\text{B}(\R^n)\simeq \R^n\times \text{SO}(n)$, thus clearly $\R^n\times
\text{Spin}(n)$ is a principal $\text{Spin}(n)$-bundle and the map given
by  $\I \times \mu:\R^n\times \text{Spin}(n) \rightarrow \R^n\times
\text{SO}(n)$  is an equivariant 2-fold cover. Thus we get a spin
structure on $\R^n$ and since $\R^n$ is contractible this is the only such
structure.

Now, if $\Gamma$ is a Bieberbach group we have a left action of
$\Gamma$ on $\text{B}(\R^n)$ given by
$\gamma\cdot(x,(w_1,\ldots,w_n))=(\gamma x,(\gamma_\ast w_1,\ldots,\gamma_\ast w_n))$. 
If $\gamma=B L_b$ then $\gamma_\ast w_j=w_j B$. Fix
$(v_1,\dots,v_n) \in \text{B}(\R^n)$. Since
$(w_1,\ldots,w_n)=(v_1g,\ldots,v_ng)$ for some $g\in
\text{SO}(n)$, we see that $\gamma_\ast w_j=(v_jg)B= v_j(B g)$,
thus the action of $\Gamma$ on $\text{B}(\R^n)$  corresponds to the
action of $\Gamma$ on $\R^n\times \text{SO}(n)$ given by
$\gamma\cdot(x,g)=(\gamma x,B g)$.

Now assume a group homomorphism is given
\begin{equation}\label{spindiagram}
\vep:\G\arr \spin \quad \text{ such that } \quad \mu \!\circ\! \vep =r,
\end{equation}
where $r(\g)=B$ if $\g=BL_b \in \G$. In this case we can lift the
left action of $\Gamma$ on $\text{B}(\R^n)$ to
$\tilde{\text{B}}(\R^n)\simeq \R^n\times \spin$ via $\gamma \cdot
(x,\tilde g)=(\g x,\vep(\g) \tilde g)$.
 Thus we have the spin structure
\begin{equation} \label{spindiag2} {\small
\begin{diagram}
  \node{\Gamma\backslash (\R^n\times \spin)}   \arrow[2]{e,t}{\overline{\I\times \mu}} \arrow{se} \node{} \node{\Gamma   \backslash (\R^n\times \son)} \arrow{sw}
  \\ \node{} \node{\Gamma\backslash \R^n} \node{}
 \end{diagram}   }
\end{equation}
for $M_\G$ since $\G \backslash \text{B}(\R^n) \simeq \text{B}(\G
\backslash \R^n)$ and $\overline{\I \times \mu}$ is equivariant.
 In this way, for each such homomorphism $\vep$, we obtain a spin
structure on $M_\G$. Furthermore, all spin structures on $M_\G$
are obtained in this manner (see \cite{Fr2}, \cite{LM}).
Throughout the paper we shall denote by $(\man,\vep)$ a spin
Bieberbach manifold endowed with the spin structure
(\ref{spindiag2}) induced by $\vep$ as in (\ref{spindiagram}).

\begin{defi} Since $r(L_\ld)=\I$ for $\ld \in \Ld$, then $\vep(\ld)=\pm1$
for any $\ld \in \Ld$. Denote by $\delta_\vep:=\vep_{|\Ld}$, the
character of $\Ld$ induced by $\vep$. We will say that a spin
structure $\vep$ on a flat manifold $M_\Gamma$ is of {\em trivial type}
if $\delta_\vep\equiv 1$. For a torus $T_\Ld$,  the only such structure is the trivial structure corresponding to $\vep \equiv 1$.
\end{defi}

\begin{rem}
The $n$-torus admits $2^n$ spin structures (\cite{Fr}). An arbitrary flat
manifold $M_\Gamma$ need not admit any.  
In \cite{MP},
we give necessary and sufficient conditions for existence in the case
when $\Gamma$  has holonomy group $\Z_2^k$ and several simple examples of flat
manifolds that can not carry a spin structure.  Also, we exhibit pairs of
manifolds, isospectral on $p$-forms for all $0\le p \le n$, where one
carries several spin (or pin$^\pm$) structures and the other admits none.
\end{rem}

\subsubsection*{Twisted spinor bundles}
Let $(L,\s)$ be the spin representation of $\text{Spin}(n)$. Then
$\dim(\s)= 2^m$, with $m=[\f n2]$. If $n$ is odd then $L$ is irreducible,
while if $n$ is even then $\s = \s^+ \oplus \s^-$, where $\s^\pm$ are
invariant irreducible subspaces of dimension $2^{m-1}$ (see Appendix).

The complex flat vector bundles over $M_\G$ are in a one to one correspondence with complex unitary representations $\rho :\G \rightarrow \text{U}(V)$. For simplicity, in this paper we shall only consider representations $\rho$ of $\G$ such that $\rho_{|\Ld}=1$.
The group $\spin$ acts on the right on
$\tilde{\text{B}}(\R^n)\times \s \otimes V$ by
$$(b,w \otimes v) \cdot \tilde g = (b\tilde g, L( \tilde{g}^{-1}) (w) \otimes v)$$ and this action defines an  equivalence relation
such that
$$((x,\tilde g), w \otimes v) \sim((x,1),L(\tilde g)(w) \otimes v)$$ for
$x\in \R^n$, $\tilde g \in \spin$, $w \otimes v\in \s\otimes V$.

There is a bundle map from the associated bundle
$\tilde{\text{B}}(\R^n)\!\times_{L\otimes\I} \s\otimes V$ onto
$\R^n\times (\s\otimes V)$, given by $\overline{((x,\tilde g),w
\otimes v)}\mapsto (x,L(\tilde g)w \otimes v)$, which is clearly
an isomorphism. Given $\rho$ as above, since
$\g\cdot\overline{((x,\tilde g),w \otimes v)} = \overline{((\g
x,\vep(\g)\tilde g),w \otimes \rho(\g)v)}$, then the corresponding
action of $\g=BL_b \in \G$ on $\R^n\times (\s \otimes V)$ is given
by
\begin{equation}\label{gamaaction}
\g\cdot (x,w \otimes v)=\big(\g x, L(\vep(\g))(w) \otimes \rho
(\g)(v)\big).
\end{equation}
 In this way we get
that the bundle $\G\backslash(\tilde{\text{B}}(\R^n)\times_{L \otimes \I}
(\s\otimes V))\arr \G\backslash \R^n=M_\G$, defined by $\vep$, is
isomorphic to
\begin{equation}   \label{spinorbundle}
S_\rho(\man,\vep):=\G \backslash(\R^n\times (\s\otimes V))\arr
\G\backslash \R^n = M_\G
\end{equation}
called the {\em spinor bundle of $\man$ with twist $V$}.

Now denote by $\G^\infty(S_{\rho}(\man,\vep))$ the space of smooth
sections of $S_{\rho}(\man,\vep)$, i.e.\@ the space of {\em spinor
fields} of $\man$. Let $\psi:\R^n\arr \R^n\times (\s\otimes V)$ be
given by $\psi(x)=(x,f(x))$ where $f:\R^n\arr\s\otimes V$ is
smooth.  We have that $\psi$ defines a section of $\G \backslash
(\R^n \times (\s \otimes V))$ if and only if, for each $\g \in
\G$, $\psi(\g x)\sim\psi(x)$, that is, if and only if $\psi(\g
x)=\tilde \g \psi(x)$, for some $\tilde \g \in \G$. Since $\G $
acts freely on $\R^n$, this is the case if and only if $\tilde
\g=\g$  and furthermore $f$ satisfies $f(\g x)=(L\!\circ\! \vep
\otimes \rho) (\g)f(x)$. In other words,
$\Gamma^\infty(S_\rho(\man,\vep))$ can be identified to the space:
 $$\{f : \R^n \rightarrow \s \otimes V \text{ smooth} \,:\,
 f(\g x)=(L\!\circ\! \vep \otimes \rho) (\g) f(x)\}.$$

In particular, in the case of $T_\Ld=\Ld\backslash\R^n$, if
$\psi(x) = (x,f(x))$ is a section of $\Ld \backslash(\R^n \times
(\s\otimes V))$ then, since we have assumed that $\rho_{|\Ld}=\I$,
in the notation of Definition 1.1,
\begin{equation} \label{equivariance}
f(x+\ld)=(L\!\circ\! \vep \otimes \rho)
(\ld)f(x)=\delta_{\vep}(\ld) f(x).
\end{equation}
 Thus $f$ is
$\delta_\vep$-equivariant. Conversely, if $f:\R^n \arr \s\otimes
V$ is $\delta_\vep$-equivariant then $\psi(x)=(x,f(x))$ is a
spinor field on $T_\Ld$.

Now $\delta_\vep \in \text{Hom}(\Ld,\{\pm1\})$, hence there exists
$u_\vep \in \R^n $ such that $\delta_\vep(\ld) = e^{2\pi i \,
u_\vep\cdot\ld}$ for all  $\ld\in \Ld$. Set
\begin{equation}\label{Lambdavep}
\Ld_\vep^\ast :=
\Ld^\ast + u_\vep
\end{equation}
  where $\Lambda^\ast$ is the dual lattice of $\Lambda$. 
If $\ld_1,\ldots,\ld_n$ is a $\Z$-basis of $\Ld$, let
$\ld_1',\ldots,\ld_n'$ be the dual basis and set
\begin{equation}\label{jepsilon}
J_\vep^\pm:=\{i \in \{1,\ldots,n\} : \vep(L_{\ld_i})=
\delta_\vep(\ld_i)=\pm1\}.
\end{equation}
We thus have,
\begin{equation}    \label{veps}
u_\vep=\tfrac{1}{2} \sum_{i\in J_\vep^-} \ld_i' \,\, \mod \Ld^*,
\end{equation}
and furthermore
\begin{equation}    \label{latis}
\Ld_\vep^\ast=\bigoplus_{j\in J_\vep^+}\Z \ld_j' \oplus\bigoplus_{j\in
J_\vep^-}(\Z+\tfrac12)\ld_j'.
\end{equation}

For  $u\in \Lambda_\vep^\ast$, $w\in \s\otimes V$ consider the function
\begin{equation} \label{fuws}
f_{u,w}(x):=f_u(x)w:=e^{2\pi i \,u\cdot x}w.
\end{equation}
It is clear that $f_{u,w}$ is $\delta_\vep$-equivariant and hence
$\psi_{u,w} (x):=(x,f_{u,w}(x))$ gives a spinor field on $T_\Lambda$.

\section{The spectrum of twisted Dirac operators.}
In this section we will introduce the Dirac operator $D_\rho$
acting on sections of the spinor bundle $S_\rho(\man,\vep)$ (see
(\ref{spinorbundle})), where $(\rho,V)$ is a finite dimensional,
unitary representation  of $\G$. We shall denote  by
$\chi_{_\rho}$ and $d_\rho $ the character and  the dimension of
$\rho$, respectively. In the main results in this section we will
give an explicit formula for the multiplicities of the eigenvalues
of twisted Dirac operators of any spin Bieberbach manifold
$(\man,\vep)$ together with a general expression for the
$\eta$-series. We shall later use this expression to compute the
$\eta$-invariant for flat manifolds with holonomy group $\Z_2^k$
(see (\ref{etainv})) and for certain $p$-dimensional flat
manifolds with holonomy group $\Z_p$, $p$ prime (see
(\ref{etainvariants})).

\subsubsection*{The spectrum of twisted Dirac operators.}
If $M_\G$ is a flat manifold endowed with a spin structure $\vep$,
let $\Delta_{s,\rho}$ denote the {\em twisted spinor Laplacian} acting on smooth sections
 of the spinor bundle $S_\rho(\man,\vep)$. That is, if $\psi(x)=(x,f(x))$
is a spinor field with $f(x)= \sum _{i=1} ^{d} f_i(x)w_i$,
$f_i:\R^n\rightarrow \C$ smooth and $\{w_i : 1\le i \le d\}$ a basis of
$\s\otimes V$ ($d=2^m d_\rho$) then
\begin{equation} \label{laplacian}
 \Delta_{s,\rho} \psi (x)= \Big( x, \sum _{i=1}^{d}
\Delta f_i (x)w_i \Big),
\end{equation}
where $\Delta$ is the Laplacian on functions on $M_\G$.

It is easy to see that for $u \in \Lambda_\vep^\ast$, $w \in\s\otimes V$,
every $f_{u,w}$ as in (\ref{fuws}) is an eigenfunction of $\Delta_{s,\rho}$ with
eigenvalue $-4\pi^2 \|u\|^2$.

\sk Let $\{e_1,\ldots,e_n\}$ be an orthonormal basis of $\R^n$ and
let $\psi$ be as above. The {\it twisted Dirac operator} $D_\rho$
is defined by
\begin{equation} \label{dirac}
D_\rho \psi (x) = \Big( x, \sum_{i=1}^n e_i  \, \f{\partial f}{\partial x_i}
\,(x) \Big)
\end{equation}
 where $e_i$ acts by $L(e_i)\otimes \I$ on $\s \otimes V$.
We will often abuse notation
and assume that $D_\rho$ acts on the function $f$ where
$\psi(x)=(x,f(x))$, writing
$D_\rho f(x) = \sum_{i=1}^n e_i \, \f{\partial}{\partial x_i} \,f(x)$.

One has that $D_\rho$ is an elliptic, first order, essentially
self-adjoint differential operator on the spinor bundle
$S_\rho(\man,\vep)$ of $M_\G$. Furthermore, $D_\rho$ has a discrete
spectrum consisting of real eigenvalues of finite multiplicity and
satisfies $D_\rho^2=-\Delta_{s,\rho}$. We will denote by
$Spec_{D}(M_\G,\vep)$
the spectrum of $D_\rho$ 
when $\rho$ is understood.
\smallskip

We now compute the action of $D_\rho$ on 
$f_{u,w}$ for  $u\in \Lambda_\vep^\ast$, $w\in \s\otimes V$. We have
$$D_\rho f_{u,w}(x)=\sum_{j=1}^n  e_j \cdot \f{\partial}{\partial x_j}
e^{2\pi i u\cdot x} w =2 \pi i\, e^{2\pi i u\cdot x} u\cdot  w =2
\pi i u\cdot f_{u,w}(x).$$ For any $u\in \R^n \smallsetminus \{0
\}$, left Clifford multiplication by $u$ on $\s$ is given by $u
\cdot w = L(u)(w)$ for $w\in \text{S}$. We fix $\langle\, ,\,
\rangle$ an inner product on $\s$ such that $L(u)$ is skew
Hermitian for every $u \in \R^n \smallsetminus \{0\}$, hence
$\langle\, ,\, \rangle$ is $\spin$-invariant. Note that $L(u)^2
=-\|u\|^2\I$. For each $u\in \Lambda_\vep^\ast$ with $\|u\|=\mu >
0$, let $\{w_j^\pm\}_{j=1}^{2^{m-1}}$ be an orthonormal basis of
the eigenspace of $L(u)$ with eigenvalue $\mp i\| u \|$. Let
$\{v_k\}_{k=1}^{d_\rho}$ be an orthonormal basis of $V$. If we set
\begin{equation} \label{fujk}
f_{u,j,k}^\pm(x) := f_{u,w_j^\pm \otimes v_k}(x) = e^{2\pi i u\cdot x}
w_j^\pm\otimes v_k,
\end{equation}
 for $u\in \Lambda_\vep^\ast$, $1\leq j \leq 2^{m-1}$
and $1\le k \le d_\rho$, then we have
\begin{equation*}
D_\rho f_{u,j,k}^\pm= \pm 2\pi \|u\|\, f_{u,j,k}^\pm,
\end{equation*}
that is, $f_{u,j,k}^\pm \in H_\mu^\pm$, the space of eigensections of $D_\rho$
with eigenvalue $\pm 2\pi\mu$.

If $u=0$, $w\ne 0$, by (\ref{equivariance}) the constant function
$f_{0,w}(x) = w $ with $w \in \s \otimes V$ defines a spinor field
if and only if
$$f_{0,w}(x+\ld) =\delta_\vep(\ld)f_{0,w}(x)$$
for
any $\ld \in \Ld$. That is, if and only if $\vep=1$, the trivial
spin structure. Moreover, $f_{0,w}$ is an eigenfunction of
$D_\rho$ with eigenvalue $0$, i.e.\@ a {\em harmonic spinor}.
Conversely, if $D_\rho f=0$ then also $\Delta_{s,\rho} f=0$ hence,
it follows from (\ref{laplacian}) that $f= f_{0,w}$, a constant
function, for some $w \in \s \otimes V$.

We shall denote by $H_0$ the space of harmonic spinors. If $n$ is even,
then  sections of the form $f_{0,w}$ with $w \in \s^+\otimes V$ (resp.\@
$\s^- \otimes V$) are often called {\em positive} (resp.\@ {\em negative})
harmonic spinors. We have $H_0=H_0^+\oplus H_0^-$, where $H_0^+$ and
$H_0^-$  respectively denote the spaces of positive and negative harmonic
spinors.

If $\mu>0$ we set $d^\pm_{\rho,\mu}(\Ld,\vep):=\text{dim}\,H_\mu^\pm$,
 the multiplicities of the eigenvalues $\pm2\pi\mu$.
If $\mu = 0$, let $d_{\rho,0}(\Ld,\vep):=\text{dim}\,H_0$ and
$d_{\rho,0}^\pm(\Ld,\vep):=\text{dim}\,H_0^\pm$, if $n$ is even.

\medskip
The next result gives  $Spec_{D}(\tor,\vep)$ for the torus $T_\Ld$ and
shows that it is determined by the cardinality of the sets
\begin{equation} \label{Ldvepmu}
\Lambda_{\vep,\mu}^\ast := \{v\in \Lambda_\vep^\ast \;:\; \|v\|=\mu\}.
\end{equation}
where $\Lambda_\vep^\ast$ is as in (\ref{Lambdavep}).

\begin{prop} \label{spectorus}
 Let $\vep$ be a spin structure on the torus $\tor \simeq \Ld \backslash
 \R^n$ and let $m=[\frac n2 ]$. In the notation above, we have:

(i) $H_0 = 
 \s \otimes V$ if $\vep = 1$ and   $H_0=0$ if $\vep\not=1$.

(ii)  If $\mu>0$, then $$H_\mu^\pm=\text{span}\{f_{u,j,k}^\pm :
u\in \Lambda_{\vep,\mu}^\ast ,\, 1\le j \le 2^{m-1},\,1\le k \le
d_\rho \}$$ with  $f_{u,j,k}^{\pm}$ as in (\ref{fujk}). The
multiplicity of the eigenvalue $\pm 2 \pi \mu$ of $D_\rho$ thus
equals
$$d_{\rho,\mu}^\pm(\Ld,\vep)= 2^{m-1}d_\rho |\Ld_{\vep,\mu}^\ast|$$
with $\Ld_{\vep,\mu}^\ast$ as in (\ref{Ldvepmu}).
Furthermore,  $\| f_{u,j,k}^{\pm}\| = \vol (T_\Lambda)^{ 1/2}$ for
each $u,j,k$.
\end{prop}

\begin{proof}
The statements in (i) are clear in light of the discussion before
the pro\-po\-si\-tion. Now let $L^2(\Ld \backslash \R^n;
\delta_\vep)$ denote the space
$$\Big\{f:\R^n \arr \C \;|\; f(x+\ld)=\delta_\vep(\ld)f(x) \text{ for } \ld \in \Ld, x \in \R^n
\text{ and } \int_{T_\Ld}|f|^2<\infty\Big\}.$$ and let $L^2(\Ld
\backslash \R^n, \s \otimes V; \delta_\vep)$ be defined similarly
by using functions with values in $\s \otimes V$.

For each $u  \in \Ld_{\vep}^\ast$  the function $f_u(x) =e^{2 \pi
i u \cdot x}$ lies in $L^2(\Ld \backslash \R^n; \delta_\vep)$. In
the case $\vep =1$, the Stone-Weierstrass theorem implies that the
set $\{f_u : u \in \Ld^ \ast \}$ is a complete orthogonal system
of  $L^2(\Ld \backslash \R^n)$. Since $\Ld_{\vep}^\ast =  \Ld^\ast
+ u_\vep$, this implies that the set $\{f_u : u \in \Ld_\vep ^\ast
\}$ is a complete orthogonal system of $L^2(\Ld \backslash \R^n;
\delta_\vep)$. Now,  if for each given $u \in \Ld_\vep^\ast$, we
choose an orthonormal basis ${\mathcal B}_u$ of $\s \otimes V$ of
eigenvectors of $L(u)\otimes \I$, then this clearly implies that
the set
\begin{equation}
\{f_{u,w}(x) : u \in \Ld_{\vep}^\ast, w \in \mathcal{B}_u \}
\end{equation}
is a complete orthogonal system of $L^2(\Ld \backslash \R^n, \s
\otimes V; \delta_\vep)$. Furthermore each $f_{u,w}$ is an
eigenfunction of $D_\rho$ with eigenvalue $\pm 2 \pi \| u\|$,
therefore (ii) follows.
\end{proof}

\begin{rem} \label{ortolattice}
 If $\Ld$ is the canonical (or cubic) lattice, then
\begin{equation} \label{weightvep}
 |\Ld_{\vep,\mu}^\ast| = |\{(m_1,\ldots,m_n)\in \Z^n :
\text{$\sum_{j\in J_\vep^+}$} m_j^2 + \text{$\sum_{j\in J_\vep^-}$}
 (m_j+\tfrac12)^2 = \mu^2\}|,
\end{equation}
where $J_\vep^\pm$ (see (\ref{jepsilon})) are computed with
respect to the canonical basis of $\R^n$. We observe that
$|J_\vep^-|$ (or $|J_\vep^+|$, see (\ref{latis})) determines the
multiplicity of the eigenvalue $\pm 2 \pi \mu$, for any $\mu
> 0$, by (\ref{weightvep}) and by the multiplicity formula in (ii) of Proposition 2.1.
The converse also holds, thus we have that
$Spec_D(T_\Ld,\vep)=Spec_D(T_\Ld,\vep')$ if and only if
$|J_\vep^-|=|J_{\vep'}^-|$.
\end{rem}

Our next goal is to obtain a formula for the multiplicities of the
eigenvalues $\pm2\pi \mu$, $\mu\ge 0$, of the Dirac operator acting on
twisted spinor fields of  an arbitrary spin flat manifold $(\man,\vep)$. 
We shall see that only the elements in $\G$ such that $B\in F_1$
(see \ref{F1}) will contribute to the formula. One of the key
ingredients in the formula will be a sign, attached to each pair
$\g,u$, with $\g=BL_b$, $u \in \Ld^\ast$ fixed by $B$. This sign
appears when comparing the conjugacy classes of two elements $x,y$
in  $\text{Spin}(n-1)$ that are conjugate in $\text{Spin}(n)$. By
Lemma \ref{conjugacy}, given two such elements $x,y$, then either
$y$ is conjugate to $x$ or to $-e_1 x e_1$ in  $\text{Spin}(n-1)$
(see Appendix). In what follows we shall write $x \sim y$ if $x,y$
are conjugate in  $\text{Spin}(n-1)$.
\smallskip

If $\g =BL_b \in \Ld\backslash \G$  and  $u \in
(\Ld_{\vep}^\ast)^B \smallsetminus \{0\}$, let $h_u \in
\text{Spin}(n)$ be such that $h_u \,u \,h_u^{-1} = \|u\| e_n$.
Therefore, since $Bu=u$, $h_u \vep(\g) h_u^{-1} \in
\text{Spin}(n-1)$, by the comments after Definition 6.5.

We shall make use of the following definitions.

\begin{defi}\label{xgama}
 Fix an element $x_\g$ in the maximal torus $T$ of $\text{Spin}(n-1)$ (see (\ref{tori})) that is conjugate in $\text{Spin}(n)$ to $\vep(\g)$.
Define  $\sigma_\vep(u,x_\g)=1$ if $h_u \vep(\g) h_u^{-1}$ is
conjugate to $x_\g$ in $\text{Spin}(n-1)$ and $\sigma_\vep(u,
x_\g)=-1$, otherwise (in this case $h_u \vep(\g) h_u^{-1}\sim -e_1
x_\g  e_1$).
\end{defi}

Note that $\sigma_\vep(u,x_\g)$ is independent of the choice of $h_u$. For simplicity, we shall simply write $\sigma(u,x_\g)$ when $\vep$ is understood.

\smallskip

\begin{rem}\label{sigmas}
If $n$ is even, then $\sigma(u,x_\g)=1$ for all $u$, by Lemma 6.2.
Also, if $\g$ is such that  $n_B>1$ then, by arguing as at the end
of the proof of Lemma \ref{conjugacy}, we see that $\sigma(u,
x_\g)=1$ for all $u$. Furthermore, if $x_\g$ is not conjugate to
$-e_1 \,x_\g \, e_1$ in $\text{Spin}(n)$, then, by the definition,
$\sigma(u, -e_1\,x_\g\, e_1)=-\sigma(u, x_\g)$. Moreover
$\sigma(-u, x_\g)=-\sigma(u, x_\g)$  and
 $\sigma(\alpha u, x_\g) = \sigma(u, x_\g)$ for any $\alpha>0$,
since we may take $h_{-u}= e_1 h_u$ and $h_{\alpha u}= h_u$,
respectively.
\end{rem}

Now, let as usual $\chi_{_\rho}$, $\chi_{_L}$ and $\chi_{_{L^\pm}}$ denote the
characters of $\rho$, $L$ and $L^\pm$, respectively. Keeping the notation above,
for each $\g = BL_b\in \G,\, \mu>0$  we set
\begin{equation}    \label{emugamas}
e_{\mu,\g,\sigma}(\delta_\vep) :=  \sum_{u \in
(\Ld_{\vep,\mu}^\ast)^B} \sigma(u,x_\g) e^{-2\pi i u\cdot b}
\end{equation}
where $(\Ld_{\vep,\mu}^\ast)^B=\{ v\in \Ld_{\vep,\mu}^\ast : Bv=v\}$.
When $\sigma=1$ we just write $e_{\mu,\g}(\delta_\vep)$
for $e_{\mu,\g,1}(\delta_\vep)$.
\smallskip

We are now in a position to prove the main result in this section.

\begin{teo}  \label{main}
Let $\G$ be a Bieberbach group with translation lattice $\Ld$ and
holonomy group $F\simeq \Ld\backslash \G$. Assume $\man=\vcp$ is a
spin compact flat manifold, with spin structure $\vep$. Then, if
$n$  is even, for each $\mu> 0$ the multiplicity of the eigenvalue
$\pm2\pi \mu$ of $D_{\rho}$ is given by
\begin{equation}    \label{multip}   d_{\rho,\mu}^\pm(\G,\vep) =
\tfrac{1}{|F|}\, \sum_{\g \in \Ld \backslash \G} \,
\chi_{_\rho}(\g)\;
e_{\mu,\g}(\delta_\vep)\;\chi_{_{L_{n-1}}}(x_\g).
\end{equation}
If $n$ is odd then
\begin{eqnarray}    \label{multipodd} \nonumber  d_{\rho,\mu}^\pm(\G,\vep)& =&
\tfrac{1}{|F|}\,\Big(\!\!\sum_{\scriptsize{\begin{array}{c}\g \in \Ld \backslash \G\\ B \not\in F_1\end{array}}}
\!\chi_{_\rho}(\g)\; e_{\mu,\g}(\delta_\vep)\;\chi_{_{L_{n-1}^\pm}}(x_\g) + \\ & &
\sum_{\scriptsize{\begin{array}{c}\g \in \Ld \backslash \G \\ B \in F_1\end{array}}}
\chi_{_\rho}(\g)\sum _{u \in (\Ld_{\vep,\mu}^\ast)^B} e^{-2\pi i u\cdot b} \;\chi_{_{L_{n-1}^{\pm \sigma(u,x_\g)}}}(x_\g)\Big)
\end{eqnarray}
with $e_{\mu,\g}(\delta_\vep)$ as in (\ref{emugamas}), $F_1$  as in (\ref{F1}), and $x_\g$, $\sigma(u,x_\g)$  as in Definition \ref{xgama}.

Let $\mu=0$. If $\vep_{|\Ld} \ne 1$ then $d_{\rho,0}(\G,\vep)=0$. If
$\vep_{|\Ld}= 1$ then
 \begin{equation}   \label{halfharmspinors}
  d_{\rho,0}(\G,\vep) = \tfrac{1}{|F|}\,
  \sum_{\g  \in \Ld \backslash   \G} \,  \chi_{_\rho}(\g)\;
  {\chi}_{_{L_{n}}}(\vep(\g))= \dim(\text{S} \otimes V)^F.
\end{equation}

We note that the summands in (\ref{multip}) and in (\ref{multipodd})  are
independent of the representative $\g$ mod $\Ld$ and of the choice of $x_\g$,  but in general the individual factors are not.
\end{teo}

\begin{proof} We proceed initially as in \cite{MR3}, \cite{MR4}.
We have that
$$ L^2(S_{\rho}(\man,\vep)) \simeq  L^2(S_{\rho}
(\tor,\vep))^\G=    \bigoplus_{\mu > 0} \Big( (H_\mu^+)^\G \oplus
(H_\mu^-)^\G \Big)\oplus {H_0}^\G.$$
Thus, 
$d_{\rho,\mu}^\pm(\G,\vep)=\dim{(H_\mu^\pm)^\G}$, for $\mu >0$,
and $d_{\rho,0}(\G,\vep)=\dim{{H_0}^\G}$.
 One has a  projection $p_{\mu}^{\pm}$ from $\hmu$ onto $(\hmu)^\Gamma$: $$p_{\mu}^{\pm}= \tfrac{1}{|\Lambda\backslash\Gamma|} \,
 \sum_{\gamma\in \Lambda\backslash\Gamma} \gamma_{|\hmu}$$
hence
$$\dim\,(\hmu)^\Gamma= \tr p_\mu^\pm = \tfrac{1}{|\Lambda\backslash\Gamma|}
  \sum_{\gamma\in \Lambda\backslash\Gamma} \tr (\gamma_{|\hmu}) $$
and similarly for  $\dim{H_0}^\Gamma$, with $H_0$ in place of $H_\mu^\pm$.

Thus, by (ii) of Proposition 2.1 we have to compute, for  $\g\in \G$,
 \begin{equation}   \label{trazagama}
 \tr \gamma_{|{\hmu}}= \f{1}{\vol(\tor)}
 \sum_{u\in \Lambda_{\vep,\mu}^\ast} \sum_{j=1}^{2^{m-1}}
  \sum_{k=1}^{d_\rho}
 \langle \gamma \cdot f_{u,j,k}^\pm, f_{u,j,k}^\pm \rangle.
 \end{equation}

  Recall, from (\ref{gamaaction}), that $\g \in \G$ acts by
  $\g\cdot \psi(x)=(\g x,(L\!\circ\! \vep \otimes \rho) (\g) f(x))$
on   $\psi(x)=(x,f(x))$. Thus there is an action of $\g$ on
$f$ given by
$$\g \cdot f(x) = (L\!\circ\! \vep \otimes \rho) (\g)
f (\g^{-1}x).$$
 Since $\g^{-1}=L_{-b}B^{-1}$ by (\ref{fujk}) we have:
\begin{equation} \label{pifjs}\begin{split}
\g \cdot f_{u,j,k}^\pm (x) & =  (L\!\circ\! \vep \otimes \rho) (\g)
\,f_{u,j,k}^\pm (\g^{-1}x)
\\ & =  e^{-2\pi i u\cdot b} \, f_{Bu}(x) \,L(\vep(\g))w_j^\pm \otimes
\rho(\g)v_k. \end{split}
\end{equation}

Let $\g_{u,j,k}^\pm  :=  \langle \g \cdot f_{u,j,k}^\pm, f_{u,j,k}^\pm
\rangle = \int_{T_\Ld} \langle \g \cdot f_{u,j,k}^\pm(x), f_{u,j,k}^\pm(x)
\rangle \;dx$. Now, using (\ref{pifjs}) we compute:
\begin{eqnarray*}
 \g_{u,j,k}^\pm & = & e^{-2\pi i u \cdot b } \int_{T_\Lambda} \langle\,
 f_{Bu}(x)L(\vep(\g))w_j^\pm \otimes
\rho(\g)v_k \,,\, f_u(x)w_j^\pm \otimes v_k \,\rangle \;dx \\& = &
e^{-2\pi i u\cdot b}  \langle\, L(\vep(\g))w_j^\pm \otimes
\rho(\g)v_k \,,\, w_j^\pm \otimes v_k  \,\rangle \,
\int_{T_\Lambda} e^{2\pi i (Bu -u)\cdot x} \;dx
\\ &= & e^{-2\pi i u\cdot b}  \langle\, L (\vep(\g)) \,w_j^\pm \otimes
\rho(\g)v_k \,,\, w_j^\pm \otimes v_k \,\rangle \; \vol(\tor) \;
\delta_{Bu,u}.
\end{eqnarray*}

In this way we get
\begin{equation*}
\tr p_\mu^{\pm}  =  \frac{1}{|F|} \sum_{\g\in \Ld\backslash \G} \!
\sum_{\scriptsize {\begin{array}{c} u\in (\Ld_{\vep,\mu}^\ast)^B
\end{array}}} \!
\sum_{k=1}^{d_\rho} \langle \rho (\g)v_k, v_k \rangle
\sum_{j=1}^{2^{m-1}} e^{-2\pi i u\cdot b} \, \langle L(\vep(\g))
w_j^{\pm}, w_j^{\pm}\rangle .
\end{equation*}

Now,  if $\gamma = BL_b\in \Gamma$ and $u \in (\Ld^*_{\vep,\mu})^B$,
then $\vep(\g) \in \text{Spin}(n-1,u)$ (see (\ref{spinu})).
Hence $L(\vep(\g))$ preserves the eigenspaces $\text{S}_u^\pm$ of $L(u)\otimes\I$ and  we
can consider the trace of $(L\!\circ\!\vep\otimes \rho) (\g)$ restricted to
$\text{S}_u^\pm \otimes V$. Thus we finally obtain
\begin{equation}            \label{premultip}
d_{\rho,\mu}^\pm(\G,\vep)  = \tf{1}{|F|}\,\sum_{\g\in \Ld\backslash \G}\,\tr
\rho(\g)\!\! \sum_{\scriptsize \begin{array}{c}   u\in
(\Ld_{\vep,\mu}^\ast)^B
\end{array}} e^{-2\pi i u\cdot b}  \;\tr
L(\vep(\g))_{|{\text{S}_u^{\pm}}}.
\end{equation}

The next task will be to compute the traces $\tr
L(\vep(\g))_{|{\text{S}_u^{\pm}}}$, showing they can be expressed
as values of characters of spin representations. The influence of
$u$ will only appear in the determination of a sign. We shall use
the element $x_\g$ and the notions introduced in Definition
\ref{xgama}.

We first note that  $\tr L(\vep (\g))_{|S_u^\pm} = \tr L(h_u\vep
(\g)h_u^{-1})_{|S_{e_n}^\pm}$. Now, we use Lemma \ref{lemaSv}
together with Definition \ref{xgama}. If $n$ odd and $h_u \vep(\g)
h_u^{-1} \sim x_\g$ then $\tr L(h_u\vep
(\g)h_u^{-1})_{|S_{e_n}^\pm} = \tr L_{n-1}^{\pm}(x_\g)$. If $h_u
\vep(\g) h_u^{-1} \not \sim x_\g$, then $h_u \vep(\g) h_u^{-1}
\sim -e_1 x_\g e_1$, hence
\begin{equation*}
\tr L(h_u\vep (\g)h_u^{-1})_{|S_{e_n}^\pm}=
\tr L_{n-1}^{\pm}(-e_1x_\g e_1)= \tr L_{n-1}^{\mp}(x_\g)
\end{equation*}
since $L(e_1)$ sends $\s^\pm$ to $\s^\mp$ orthogonally. For $n$
even we proceed similarly, using (\ref{spinrestrictions}). Thus we
obtain:
\begin{equation}\label{tracesu}
\tr L(\vep (\g))_{|S_u^\pm} 
= \left\{\begin{array}{ll}\tr  L_{n-1}(x_\g)& n \text{ even} \sk \\
 \tr L_{n-1}^{\pm\sigma(u, x_\g)}(x_\g)& n \text{ odd}.\end{array} \right.
\end{equation}

Substituting (\ref{tracesu}) in (\ref{premultip}) we get that
$d_{\rho,\mu}^\pm(\G,\vep)$ equals
\begin{equation*}
\tf{1}{|F|}\,\sum_{\g\in \Ld\backslash \G}\,\chi_\rho(\g) \!\!
\sum_{\scriptsize \begin{array}{c}   u\in (\Ld_{\vep,\mu}^\ast)^B
\end{array}} e^{-2\pi i u\cdot b}
\left\{\begin{array}{ll} \! \tr  L_{n-1}(x_\g)&  \quad n \text{ even} \sk \\
 \! \tr L_{n-1}^{\pm \sigma(u,x_\g)}(x_\g)&   \quad n \text{ odd} \end{array} \right.
\end{equation*}
as asserted in formula (\ref{multip}) for $n$ even. If $n$ is odd,
then, by separating the contributions of the elements $\g =BL_b$
with $B \in F_1$ from those with $B \notin F_1$ (see Remark
\ref{sigmas}), we arrive at formula (\ref{multipodd}).

\sk In the case when $\mu=0$ we may proceed in a similar  way. If
we identify $w\otimes v \in \s\otimes V$ with the constant
function $f_{0, w\otimes v}$, then for $\g \in \G$ we have
\begin{equation*}
\text{tr} \, \g_{|{H_0}} = \tfrac 1 {\text{vol}\, T_\Ld} \sum _{j=1}^{2^m}
\sum_{k=1}^{d_\rho} \langle \g \cdot w_j \otimes v_k, w_j\otimes v_k \rangle =
\chi_{_\rho}(\g) \; \chi_{_{L_n}}(\vep (\g)).
\end{equation*}
Thus
\begin{equation*}
\text{dim}\,(H_0^\G) =\tfrac 1 {|F|} \sum_{\g \in \Ld \backslash \G}
\text{tr}\,\g_{|{H_0}}=\tfrac 1 {|F|}\sum_{\g \in \Ld \backslash \G}
\chi_{_\rho}(\g) \; \chi_{_{L_n}}(\vep (\g))
\end{equation*}
as claimed. Concerning the last equality in the theorem,  we know
that
$$H_0=\{ f_{0,w\otimes v} : L\!\circ\! \vep \otimes \rho (\g)w\otimes v = w \otimes v \} \simeq
(\s \otimes V)^F,$$
since $F\simeq \Ld \backslash \G$ and $\Ld$ acts trivially for $\vep$ of
trivial type.
\end{proof}

\begin{coro} \label{simetspec}
Let $(M_\G,\vep)$ be a spin compact flat manifold of dimension
$n$. If $n$ is even, or if $n$ odd and $n_B>1$ for every $\g =
BL_b \in \G$ (i.e. $F_1=\emptyset$), then the spectrum of the
twisted Dirac operator $D_\rho$ is symmetric.
\end{coro}
\begin{proof}
If $n$ is even, the assertion is automatic by (\ref{multip}). If
$n$ is odd, by (\ref{multipodd}) and (\ref{spinpmchars}) we have
that $d_{\rho,\mu}^+(\G,\vep) - d_{\rho,\mu}^-(\G,\vep)$ equals
\begin{equation}\label{difference}
 \tf{(2i)^{m}}{|F|}\,\sum_{\g\in \Ld\backslash \G}\,\chi_\rho(\g)\,
e_{\mu,\g}(\delta_\varepsilon) \prod_{j=1}^m  \sin t_j(x_\g)
\end{equation}
where $x_\g=x(t_1(x_\g),\dots, t_m(x_\g))$, in the notation of (\ref{torielements}).
It is clear that $\g$ satisfies $n_B>1$ if and only if $t_k(x_\g) \in \pi \Z$, for some $k$, hence in this case  $\prod_{j=1}^m  \sin t_j(x_\g)=0$.
Thus, if $n_B>1$ for every $\g$, (\ref{difference}) implies that the Dirac spectrum is symmetric.
\end{proof}
\begin{rem}
In specific examples, the expressions for the multiplicities of eigenvalues in the theorem, can be made more explicit by  substituting $\chi_{_{L_{n-1}}}(x_\g)$, $\chi_{_{L^\pm_{n-1}}}(x_\g)$ or $\chi_{_{L_{n}}}(\vep(\g))$ by the values
given in Lemma \ref{spincharacters} in terms of products of cosines or sines of the  $t_j(x_\g)$, where $x_\g =x(t_1(x_\g),\dots, t_m(x_\g))\!$ as above.
\end{rem}
\begin{rem} \label{spinlapmultip}
Set $d_{\rho,\mu}(\G,\vep):=d_{\rho,\mu}^+(\G,\vep) +
d_{\rho,\mu}^-(\G,\vep)$ for $\mu>0$. Note that, since $D_\rho^2=-\Delta_{s,\rho}$, then
$d_{\rho,\mu}(\G,\vep)$ is just the multiplicity of the eigenvalue
$4\pi^2\mu^2$ of $-\Delta_{s,\rho}$.
 We thus have:
\begin{equation}    \label{spinlapmultipform}
d_{\rho,\mu}(\G,\vep) =   \tfrac{\alpha(n)}{|F|}\, \sum_{\g\in
\Ld\backslash \G}\,  \chi_{_\rho}(\g)\;
e_{\mu,\g}(\delta_\vep) \; \chi_{_{L_{n-1}}}(\vep(\g))
\end{equation}
where $\alpha(n)=1$ or $2$ depending on whether $n$ is odd or even.
Clearly, (twisted) Dirac isospectrality implies (twisted) spinor Laplacian isospectrality, but we shall see the converse is not true (see Example \ref{example2}).
\end{rem}

\subsubsection*{Spectral asymmetry and $\eta$-series} We decompose
 $Spec_D(M_\G,\vep)=\mathcal{S} \,\dot\cup\, \mathcal{A}$ where
$\mathcal{S}$ and $\mathcal{A}$ are the symmetric and the
asymmetric components of the spectrum, respectively. That is, if
$\ld=2\pi\mu$, $\ld\in\mathcal{S}$ if and only if
$d_{\rho,\mu}^+(\G,\vep)=d_{\rho,\mu}^-(\G,\vep)$. We say that
$Spec_D(M_\G,\vep)$ is symmetric if $\mathcal{A}=\emptyset$. In
this case, the positive spectrum $Spec_D^+(M_\G) = \{ \ld \in
Spec_D(\man) : \ld>0 \}$ and $H_0$ determine the whole spectrum
$Spec_D(M_\G)$. The symmetry of the spectrum of the Dirac operator
depends on $\G$ and also on the spin structure $\vep$ on $M_\G$.

Our next goal is to derive an expression for the $\eta$-series
$\eta_{(\G,\rho,\vep)}(s)$, for a general Bieberbach manifold $M_\G$ with a
spin structure $\vep$. Consider
\begin{equation}    \label{etaseries}
\sum_{\scriptsize \begin{array}{c} \ld \in
Spec_{D}(M_\G,\vep) \\ \ld \not= 0 \end{array}} \frac{\text{sgn}(\ld)}{|\ld|^{s}}
= \frac{1}{(2\pi)^s} \sum_{\mu \in \frac{1}{2\pi}\mathcal{A}}
\frac{d_{\rho,\mu}^+(\G,\vep) - d_{\rho,\mu}^-(\G,\vep)}{|\mu|^s}.
\end{equation}
It is known that this series 
converges absolutely for $Re(s)>n$ and defines a holomorphic
function  $\eta_{(\G,\rho,\vep)}(s)$ in this region, having a
meromorphic continuation to $\C$, with simple poles (possibly) at
$z=n-k$, $k\in \N_0$, that is holomorphic at $s=0$ (\cite{APS} and
\cite{Gi}). One defines the eta-invariant of
 $M_\G$ by $\eta_{(\G,\rho,\vep)}:= \eta_{(\G,\rho,\vep)}(0)$.
   By Corollary \ref{simetspec}, $\eta(s) \equiv 0$ if $n=2m$ or if $n=2m+1$ and $F_1=\emptyset$.
Furthermore, it is known that if $n\not\equiv 3\mod(4)$ then
$\eta(s)\equiv 0$ for every Riemannian manifold $M$. In what
follows we will often write $\eta(s)$ and $\eta$ in place of
$\eta_{(\G,\rho,\vep)}(s)$ and $\eta_{(\G,\rho,\vep)}$, for
simplicity.

\begin{prop}
Let $\G$ be a Bieberbach group of dimension $n=4r+3$ (thus $m=2r+1$)
with holonomy group $F$ and let $\vep$ be a
 spin structure on  $M_\G=\vcp$.
Then the $\eta$-series    of $M_\G$ is given by
 \begin{equation}   \label{etaodd}
  \eta_{(\G,\rho, \vep)}(s) = \tfrac{(2i)^m}{|F|(2\pi)^s}
  \sum_{\scriptsize{\begin{array}{c} \g \in \Ld\backslash \G \\ B\in F_1 \end{array}}}  \chi_{_\rho}(\g) \; \Big(\prod_{j=1}^m \sin t_j(x_\g)\Big)
 \sum_{\mu \in \frac{1}{2\pi}\mathcal{A}}
  \frac{e_{\mu,\g,\sigma}(\delta_\vep)}{|\mu|^s}
 \end{equation}
where  $e_{\mu,\g,\sigma}(\delta_\vep) =  \sum_{u \in
(\Ld_{\vep,\mu}^\ast)^B} \sigma(u,x_\g) e^{-2\pi i u\cdot b} $  as
in (\ref{emugamas}) and, in the notation of (\ref{torielements}),
$x_\g=x(t_1(x_\g),\dots, t_m(x_\g))$.
\end{prop}
\begin{proof}
By (\ref{multipodd}) and Corollary 2.6, we get that
$d_{\rho,\mu}^+ (\G,\vep) - d_{\rho, \mu}^- (\G,\vep)$ equals
\begin{eqnarray*} \label{difference2}
 & &\tfrac{1}{|F|} \sum_{\scriptsize{\begin{array}{c} \g \in \Ld\backslash \G \\ B\in F_1 \end{array}}}
 \chi_{_\rho}(\g) \;
\sum_{u\in (\Ld_{\vep,\mu}^\ast)^B} e^{- 2\pi i u\cdot b}\; \big(
\chi_{L_{n-1}^{\sigma(u,x_\g)}} - \chi_{L_{n-1}^{-\sigma(u,x_\g)}} \big) (x_\g)\\
&=& \tfrac{1}{|F|} \sum_{\scriptsize{\begin{array}{c} \g \in
\Ld\backslash \G \\ B\in F_1 \end{array}}} \chi_{_\rho}(\g) \;
\sum_{u\in (\Ld_{\vep,\mu}^\ast)^B} e^{- 2\pi i u\cdot
b}\sigma(u,x_\g)\; \big( \chi_{L_{n-1}^+} - \chi_{L_{n-1}^-} \big)
(x_\g).
\end{eqnarray*}
Now, using (\ref{spinpmchars}) we get
\begin{equation*}
d_{\rho,\mu}^+ (\G,\vep) - d_{\rho,\mu}^- (\G,\vep) =
\tfrac{(2i)^m}{|F|} \sum_{\scriptsize{\begin{array}{c} \g \in
\Ld\backslash \G \\ B\in F_1 \end{array}}} \chi_{_\rho}(\g) \;
e_{\mu,\g,\sigma}(\delta_\vep)\;  \prod_{j=1}^m \sin t_j(x_\g).
\end{equation*}
By substituting this last expression in (\ref{etaseries}) the proposition follows.
\end{proof}

\begin{rem}
(i) In Proposition 3.4 we will give a very explicit expression for the eta
series and  will compute the eta invariant for a general flat manifold
with holonomy group $\Z_2^k$. We shall see that we may have $\eta=0$ or
$\eta\ne 0$, depending on the spin structure.

(ii) Let $g\in \spin$. If $n$ even, $\chi_{L_n^+}(g) -
\chi_{L_n^-}(g)=\text{Str}\, L_n(g)$, is the {\em supertrace} of
$L_n(g)$. Furthermore, by Proposition 3.23 in \cite{BGV},
one has the expression:
\begin{equation*} \label{supertrace}
 \text{Str}\, L_n(g) = i^{-n/2}\; \text{sgn}(g) \;
|\!\det(\I_{n-1} - \mu(g))|^{1/2}
\end{equation*}
where $\text{sgn}(g)\in \{\pm 1\}$ is defined in \cite{BGV} and
$\mu$ is the covering map (\ref{mucovering}).

(iii) Note that some authors use  $\tf 12(\eta(0) + d_0)$ instead
of $\eta(0)$ as the definition of the  $\eta$-invariant of
$(M,\vep)$.
\end{rem}


\section{The twisted Dirac spectrum of $\Z_2^k$-manifolds.}
In this section we shall look at the case of $\Z_2^k$-manifolds, a
very rich class of flat manifolds. In this case, the holonomy
group $F\simeq \Z_2^k$, but the holonomy action need not be
diagonal in general. Already in the case when $k=2$ it is known
that there are infinitely many indecomposable holonomy actions and
infinitely many of these give rise to (at least) one Bieberbach
group (see \cite{BGR}). Giving a classification for $k\ge 3$ is a
problem {\em of wild type}. The case $k=n-1$ corresponds to the so
called generalized Hantzsche-Wendt manifolds, studied in
\cite{MR}, \cite{RS} and \cite{MPR}. This class is still very rich
and will be used in Section 4 to construct large families of Dirac
isospectral manifolds pairwise non \! homeomorphic to each
other\!.

If $\G$ has holonomy group $\Z_2^k$, then
$\G=\langle \g_1,\dots,\g_k, \Ld \rangle$  where
$\g_i=B_iL_{b_i}$, $B_i \in \on$, $b_i \in \R^n$, $B_i\Ld=\Ld$,
$B_i^2=\I$ and $B_i B_j=B_j B_i$, for each $1\leq i,j\leq k$. We
shall see that, somewhat surprisingly,  for these manifolds the
multiplicity formulas in Theorem \ref{main} take  extremely simple
forms (see (\ref{multeven}) and (\ref{multodd})). This will allow
to exhibit, in the next section, large Dirac isospectral sets.
Also, we will be able to characterize all  manifolds having
asymmetric Dirac spectrum and to obtain very explicit expressions
for the $\eta$-series and the $\eta$-invariant.
\sk

Let $F_1$ be as in (\ref{F1}). In the case of $\Z_2^k$-manifolds,
$F_1$ is the set of $B\in \text{O}(n)\cap r(\G)$ such that $B$ is
conjugate in $\text{O}(n)$ to  the diagonal matrix
$diag(-1,\dots,-1,1)$. The next lemma will be very useful in the
proof of Theorem \ref{main2}, the main result in this section.

\begin{lema}                \label{diagB}
If $\G$ is a Bieberbach group with translation lattice $\Lambda$ and
holonomy group $\Z_2^k$ then the elements in $F_1$ can be simultaneously
diagonalized in $\Lambda$, that is, there is a basis $f_1,\dots, f_n$ of
$\Lambda$ such that $B f_j = \pm f_j, \, 1\le j \le n$,  for any $BL_b \in
\G$ with $n_B=1$. Furthermore $\G$ can be conjugated by some $L_\mu, \,
\mu \in \R^n,$ to a group $\G'$ such that $2b\in \Ld$ for any $BL_b \in
\G'$ with $n_B=1$.
\end{lema}

\begin{proof}
Let $BL_b \in \G$ with $n_B=1$. Then $B\Ld = \Ld$ and it is a well
known fact that $\Ld$ decomposes, with respect to the action of
$B$, as $\Ld = \Ld_1 \oplus \Ld_2$ where $\Ld_1$ (resp.\@ $\Ld_2$)
is a direct sum of integral subrepresentations of rank 1 (resp.\@
2). Here   $B$ acts diagonally on $\Ld_1$, whereas $\Ld_2$ is a
direct sum of  subgroups on which  $B$
acts by $J=\left[ \begin{smallmatrix} 0 & 1 \\
1 & 0  \end{smallmatrix} \right]$. Now, it is not hard to check
that since $\G$ is Bieberbach, then the orthogonal projection
$2b^+ \in \Lambda$ of $2b$ onto $(\R^n)^B$ can not lie in $\Ld_2$,
otherwise some element of the form $BL_{b+\lambda}$, with $\lambda
\in \Ld$, would be of finite order (see \cite{DM}, Proposition
2.1). Thus,  the multiplicity of the eigenvalue 1 for $B$ is at
least 1 on $\Ld_1$. If furthermore  $\Ld_2\ne 0$, then it would be
at least 1 on $\Ld_2$, hence  $n_B \ge 2$, which is not possible
since $B \in F_1$. Thus  $\Ld_2=0$ for every $B\in F_1$, therefore
each such $B$ can be diagonalized in a basis of $\Ld$ (see Lemma
3.3 in \cite{RS} for a different proof).

Our next task is to show that this can be done simultaneously in $\Ld$ for
all elements in $F_1$. To show this, we enumerate  the
elements in $F_1$: $B_1,\dots, B_r$. Let $f_1,\dots, f_n$ be a basis of $\Ld$ diagonalizing
$B_r$. After reordering we may assume that $B_r f_n=f_n$ and $B_r f_j =
-f_j$ for $1 \le j \le n-1$. Clearly $f_n$ is orthogonal to $f_j$ for $
j<n$. Now, for $i=1, \dots, r-1$, $B_i f_n =\pm f_n$ since $B_i$ commutes
with $B_r$, and actually $B_i f_n =- f_n$, since otherwise we would have
$B_i =B_r$, since $B_i \in F_1$. Now, $B_1, \dots, B_{r-1} \in F_1$ and
leave $\Ld':=\Z\text{-span}\{f_1,\dots,f_{n-1}\}$ invariant. Hence by
induction they can be simultaneously diagonalized in some $\Z$-basis of
$\Ld'$. Putting  this basis together with $f_n$ we get a basis of $\Ld$
that diagonalizes $B_1, \dots, B_r$. The proof of the second assertion now
follows in the same way as that of Lemma 1.4 in \cite{MR4}.
\end{proof}

We are now in a position to prove the main result in this section.
We shall see that for a spin $\Z_2^k$-manifold, only the identity element, $\I$, and possibly one element $BL_b \in \G$ with $n_B = 1$ can give a nonzero contribution to the multiplicity formulas.

We note that given  $\g=BL_b\in \Gamma \smallsetminus \Ld$, since
$B^2=\I$, then $B$ is conjugate in $\text{SO}(n)$ to a diagonal matrix of
the form $\text{diag}(\underbrace {-1,\dots,-1}_{2h}, 1,\dots,1)$ with
$1\le h\le m$ where $m=[\frac n2]$, hence $\vep(\gamma)$ is conjugate in
$\text{Spin}(n)$ to $\pm g_h$ for some $1 \le h\le
m$, where
\begin{equation}\label{g_h}
 g_h:= e_1 \cdots e_{2h} \in \text{Spin}(n).
\end{equation}
 Moreover, if $n=2m$ then $h<m$, since $B=-\I$ cannot occur for $BL_b
\in \G$, $\Gamma$ being torsion-free. Thus, it follows that $g_h
\in \text{Spin}(n-1)$ and since $g_h$ and $-g_h$ are conjugate
(see Corollary \ref{tracespin}), then we may take $x_\g$ in
Definition \ref{xgama} to be equal to $g_h$, with $h$ depending on
$\g$. From now on in this section {\it we shall thus \it assume}
that $x_\g =g_h$.

We will need the fact that (see Lemma \ref{spincharacters}), if $n=2m$, then
\begin{equation}\label{spinghchars}
\chi_{_{L^\pm_n}} (g_h)  = \left \{\begin{array}{ll}\pm 2^{m-1}i^m
& \quad h=m \sk \\ 0 & \quad 1\le h <m
\end{array}\right.
\end{equation}
and furthermore, for $n$ even or odd one has
\begin{equation} \label{spinghchartot}
\chi_{_{L_n}} (g_h)  = 0 \qquad \text{ for any } 1\le h \le m.
\end{equation}

\begin{teo}     \label{main2}
Let $(\man,\vep)$ be a spin $\Z_2^k$-manifold of dimension $n$,
and let $F_1$ be as in (\ref{F1}).
\smallskip

(i) If $F_1=\emptyset$, then $Spec_{D_\rho}(\man,\vep)$ is symmetric and the
multiplicities of the eigenvalues $\pm2\pi\mu$, $\mu>0$, of $D_\rho$ are given by:
\begin{equation}    \label{multeven}
d_{\rho,\mu}^\pm(\G,\vep) = 2^{m-k-1} \,d_\rho\,  |\Ld^\ast_{\vep,\mu}|.
\end{equation}

(ii) If $F_1 \ne \emptyset$ then $Spec_{D_\rho}(\man,\vep)$ is
asymmetric if and only if the following conditions hold: $n=4r+3$
and there exists $\g=BL_b$, with $n_B=1$ and
$\chi_{_\rho}(\g)\not=0$, such that $B_{|\Ld} =-\delta_\vep \I$.
In this case,  the asymmetric spectrum is the set of eigenvalues
\begin{equation*}
{\mathcal A} = Spec_{D_\rho}^A(\man,\vep)= \{\pm 2 \pi \mu_j :
\mu_j =(j+\tfrac12){\|f\|}^{-1},\, j \in \N_0\}
\end{equation*}
where $\Ld^B= \Z f$ and if we put $\sigma_\g:=\sigma((f\cdot 2b)f,g_m)$ we have:
\begin{equation}\label{multodd}
d_{\rho,\mu}^\pm(\G,\vep)  = \left\{\!\! \begin{array}{lc} 2^{m-k-1}\,
\big( d_\rho  \, |\Ld^\ast_{\vep,\mu}| \pm 2 \sigma_\g (-1)^{r+j}
 \chi_{_\rho}(\g) \big) & \mu=\mu_j,
\sk \\ 2^{m-k-1}\, d_\rho \, |\Ld^\ast_{\vep,\mu}| & \mu \ne \mu_j
\end{array} \right.
\end{equation}

If $Spec_{D_\rho}(M_\G,\vep)$ is symmetric then $d_{\rho,\mu}^\pm(\G,\vep)$ is given by (\ref{multeven}).

(iii) The number of independent harmonic spinors is given by
$$d_{\rho,0}(\G,\vep)=\left\{ \begin{array}{ll} 2^{m-k} d_\rho & \qquad
\text{if $\vep_{|\Ld}=1$} 
 \sk \\  0 & \qquad \text{otherwise.}
\end{array} \right.$$
If $k>m$ then $\man$ has no spin structures of
trivial type, hence $M_\G$ has no harmonic spinors. Furthermore, if $\man$
has exactly $2^m d_\rho$ harmonic spinors then $\man = \tor$ and $\vep=1$.
\end{teo}

\begin{proof}
We first note that the contribution of $\I \in F$ to the multiplicity
formulas (\ref{multip}) and (\ref{multipodd}) is given by $2^{m-k-1}  d_\rho\,
e_{\mu,\I}(\delta_\vep)= 2^{m-k-1} d_\rho \,
|\Ld^\ast_{\vep,\mu}|$. Hence, when no element in $F$ other than $\I$
gives a nonzero contribution then (\ref{multeven}) holds.

If $F_1 = \emptyset$, then,  for any $\g=BL_b\in \G\smallsetminus \Ld$, we
have $\vep(\g)\sim g_h$ with  $h<m$. Thus, 
$\chi_{_{L_n}}(\vep(\g))=0$ if $n$ odd and $\chi_{_{L_n}}^\pm(\vep(\g))=0$ if $n$ even,  for any
$\g\in \G\smallsetminus\Ld$. Thus, in this case, only $\I \in F$
contributes to (\ref{multip}) and (\ref{multipodd}) and hence (i) follows. This implies that, other than $\I$, only the elements in $F_1$ can give a nonzero
contribution to the multiplicities and furthermore, this can happen only if
$n$ is odd.

Now, assume that $n=2m+1$ and $F_1 \ne \emptyset$. If
$\g=BL_{b}\in \G\smallsetminus \Ld$, with $n_B=1$ (hence
$x_\g=g_m$), we know by Lemma \ref{diagB} that there exists a
basis $f_1,\dots,f_n$ of $\Lambda$ such that $B$ is diagonal in
this basis. After reordering the basis elements we may assume that
\begin{equation}\label{Bconditions}
B f_j= -f_j \text{ for } 1\le j <n,\quad B f_n= f_n, \quad b
\equiv \frac 12 f_n \mod \Ld.
\end{equation}
Let $f'_1,\dots,f'_n \in \Lambda^\ast$ be the dual basis of
$f_1,\dots,f_n$. It is clear that also $B f'_j= - f'_j$ for $1\le
j <n$ and  $f'_n = \frac{f_n}{\|f_n\|^2}$, thus $B f'_n= f'_n$.

Let as usual $\Ld^\ast_\vep=  \Lambda^\ast + u_\vep$, with $u_\vep
= \sum_j c_j f'_j$ and $c_j \in \{ 0, \frac 12\}$ for each $j$. If
the contribution of $B$ to (\ref{multipodd}) is non-trivial, then
$(\Ld^\ast_\vep)^{B}\not=\emptyset$. Thus, there exists
 $u= \ld' + u_\vep $ with $\ld'=\sum_j d_j f'_j, \, d_j \in \Z$
 and such that $Bu = u$.
 This says
that for $1\le j \le n-1$, we have $c_j +d_j = -c_j -d_j$ and hence $c_j
=0$ for $1\le j \le n-1$. On the other hand, $c_n $ equals 0 or $\frac
12$, that is, $u_\vep =0$ or $u_\vep = \frac 12 f'_n$.

If $u_\vep =0$ we have $\Ld^\ast_\vep=\Ld^\ast$ and thus $e^{2 \pi i
u\cdot b} = e^{-2 \pi i u \cdot b}$ for any $u \in \Ld^\ast$, since $b \in \tf 12 \Ld$, by  Lemma
\ref{diagB}. 
We claim that $\sum _{u} e^{-2\pi i u\cdot b} \;\chi_{_{L_{n-1}^{\pm \sigma(u,x_\g)}}}(x_\g)=0$ where the sum is taken over
$u \in (\Ld_{\vep,\mu}^\ast)^B$.
Indeed, putting together the contributions of $u$ and $-u$ in the expression
above, by Remark \ref{sigmas} and (\ref{spinghchartot}) we get
\begin{equation*}\label{signos1}
e^{-2 \pi i u\cdot b}\,\Big(\chi_{_{L_{n-1}^{\pm \sigma(u,g_m)}}}
+ \chi_{_{L_{n-1}^{\pm \sigma(-u,g_m)}}}\Big) (g_m)= e^{-2 \pi i
u\cdot b}\,\sigma(u,g_m)\, \chi_{_{L_{n-1}}}(g_m)=0.
\end{equation*}
Hence,  we conclude that if $u_\vep=0$, the
contribution of  $\g =BL_b\in \G\smallsetminus \Ld$ to (\ref{multipodd}) is
zero.
\smallskip

Now consider the case  $u_\vep=\frac12 f_n'$, that is,
 $(\Ld^\ast_\vep)^B= (\Z + \frac 12)f'_n$. Hence, since
$\delta_\vep (\ld)=e^{2\pi i u_\vep \cdot \ld}$, then
\begin{equation}\label{goodvep} \vep(L_{f_j}) =1 \quad
(1\le j \le n-1), \qquad \mbox{ and } \qquad \vep(L_{f_n}) =-1.
\end{equation}
  Furthermore we note that, since
$$-1=\vep(L_{f_n}) = \vep(\g^2) = \vep(\g)^2= (\pm e_1\dots
e_{2m})^2= (-1)^m$$
 it  follows that $m$ is necessarily odd, hence $m=2r+1$ and $n=4r+3$.

Now (\ref{Bconditions}) and (\ref{goodvep}) say that only if
$B_{|\Ld} =-\delta_\vep \I$ can $\g=BL_b$ give a nonzero
contribution.

Furthermore, since $(\Ld^\ast_{\vep})^{B}= (\Z + \frac 12)f'_n$,
then, for fixed $\mu >0$, $(\Ld^\ast_{\vep,\mu})^{B}\ne \emptyset$ if and
only if $\mu = \mu_j := (j + \frac 12)\|f_n\|^{-1}$ with $j \in \N_0$. For
$\mu=\mu_j$, we have that  $(\Ld^\ast_{\vep,\mu_j})^{B}=\{\pm u_j\}$ where
$u_j=(j+\frac 12)f'_n$ with $j= \| f_n\|\mu_j -\frac 12 \in \N_0$.
Again, putting together the contributions of $u_j$ and $-u_j$, we
get that the sum over $(\Ld^\ast_{\vep,\mu})^{B}$ in
(\ref{multipodd}) equals
\begin{eqnarray*}
& &e^{-2 \pi i u_j\cdot b}\,\chi_{_{L_{n-1}^{\pm \sigma(u_j,g_m)}}}(g_m) +e^{2 \pi i u_j\cdot b}\, \chi_{_{L_{n-1}^{\pm \sigma(-u_j,g_m)}}} (g_m) \\
&= &\big(e^{-2 \pi i u_j\cdot b}- e^{2 \pi i u_j\cdot b}\big)\,\sigma(u_j,g_m) \, \chi_{_{L_{n-1}^\pm}}(g_m)\\
&= &\mp 2^m i^{m+1}\,\sigma(u_j,g_m) \sin(2\pi u_j\cdot b).
\end{eqnarray*}
where we have used that $\chi_{_{L_{n-1}^-}}(g_h) =
-\chi_{_{L_{n-1}^+}}(g_h)$ by (\ref{spinghchars}). If $\Ld^B=\Z
f$, then one has that $f =\pm f_n$.
Now one verifies that  $\sigma(u_j,g_m)= \sigma((f\cdot 2b)f,g_m)= \sigma_\g$. 
Hence,
$$\sigma(u_j,g_m) \sin(2\pi u_j\cdot b)=  \sigma_\g \sin(\pi (j+\tfrac 12)) =
\sigma_\g \,(-1)^{j}$$ since $b \equiv \frac 12 f_n \mod \Ld$.

Since $m=2r+1$, we finally get
that the contribution of $\g$ to the multiplicity of the eigenvalue
$\pm 2\pi \mu_j$ is given by
\begin{equation*}
\pm 2^{m-k}\sigma_\g \,(-1)^{r+j} \chi_{_\rho}(\g)
\end{equation*}

The above shows that if an element $\g'=B'L_{b'} \in F_1$ gives a
nonzero contribution, then $B'_{|\Ld} =-\delta_\vep \I$,
 hence $B'=B$, and $b'\equiv\frac 12 f_n, \mod \Ld$. Since only $\I$ and $\g
=BL_b$ give a contribution to the multiplicity formula, this
completes the proof of (ii).

Finally, the first assertion in (iii) follows immediately from
(\ref{halfharmspinors}) and the remaining assertions are direct
consequences of the first.
\end{proof}

\begin{rem} \label{remark}
(i) Except for the very special case described in (ii) of the theorem,
 the twisted Dirac spectrum of $\Z_2^k$-manifolds is symmetric
  and the multiplicities are given by the simple formula  (\ref{multeven}).
In this case,  the multiplicities of Dirac eigenvalues for $\man$ with a
spin structure $\vep$  are determined by the multiplicities for the covering
torus $T_\Lambda$ with the restricted spin structure $\vep_{|\Ld}$. Indeed, we have
$$d_{\rho,\mu}^\pm(\G, \vep)= 2^{-k} \,d_{\rho,\mu}^\pm(\Ld, \vep_{|\Ld}).$$

(ii) Note that, for each  fixed $\rho$, all spin
$\Z_2^k$-manifolds $(\man,\vep)$ having asymmetric Dirac spectrum
and having  the same covering torus $(\tor,\vep_{|\Ld})$ are
$D_\rho$-isospectral to each other.
\end{rem}

As an application of Theorem \ref{main2} we shall now compute the
$\eta$-series and the $\eta$-invariant for any $\Z_2^k$-manifold.

\begin{prop}            \label{prop.eta2k}
Let $(M_\G,\vep)$ be a spin $\Z_2^k$-manifold of odd dimension \linebreak $n=4r+3$ (thus $m=2r+1$).
If $\text{Spec}_{D_\rho}(M_\G,\vep)$ is asymmetric then, in the notation of Theorem \ref{main2}, we have:
\begin{equation} \label{etaz2kmanif}
\eta_{(\G,\rho,\vep)}(s) =  (-1)^{r}\,  \sigma_\g \,\chi_{_\rho}(\g) \,2^{m-k+1}\,\frac{\|   f \|^s}{(4\pi)^s}  \big( \zeta(s, \tfrac 14) - \zeta(s,
\tfrac34)\big)
\end{equation}
where $\zeta(s,\alpha)=\sum_{j=0}^\infty \tf{1}{(j+\alpha)^s}$
denotes the Riemann-Hurwitz zeta function for $\alpha\in (0,1]$
and $\sigma_\g \in \{\pm 1\}$ is as defined in $(ii)$ of Theorem
\ref{main2}.

Therefore, $\eta_{(\G,\rho,\vep)}(s)$ has an  analytic continuation to
$\C$ that is everywhere holomorphic. Furthermore,
\begin{eqnarray}    \label{etainv}
  \eta_{(\G,\rho,\vep)}(0)  & = &     (-1)^{r} \sigma_\g \, \chi_{_\rho}(\g)\, 2^{m-k}, \\
\label{p.etadif} \eta'_{(\G,\rho,\vep)}(0) & = & (4 \log
\G(\tfrac14) + \log\| f\| - 3\log(2\pi) ) \, \eta(0).
\end{eqnarray}
\end{prop}

\begin{proof}
We shall use Theorem \ref{main2}  in the case of the special
spin structure when the spectrum is not symmetric, otherwise,
$\eta(s)= 0$. We have
\begin{equation*}
\eta_{(\G,\rho,\vep)}(s) =  \frac 1{(2\pi)^s}\sum _{j=0}^\infty
\frac{d_{\rho,\mu_j}^+(\G,\vep)-d^-_{\rho,\mu_j}(\G,\vep)}{|\mu_j|^s}.
\end{equation*}
Now from formula (\ref{multodd}) we have that
\begin{equation*}
d^+_{\rho,\mu_j}(\G,\vep)-d^-_{\rho,\mu_j}(\G,\vep)=
(-1)^{r+j} \, \sigma_\g \, \chi_{_\rho}(\g) \,2^{m-k+1}.
\end{equation*}

Thus, if $Re(s) > n$
\begin{eqnarray*}
 \eta_{(\G,\rho,\vep)}(s)& = & (-1)^{r} \sigma_\g
 \,\chi_{_\rho}(\g)\,
  2^{m-k+1} \,\frac{||f||^s}{(2\pi)^s}  \sum_{j=0}^\infty \frac {(-1)^j}{(j+\frac 12)^s}
 \\ & = & (-1)^{r} \sigma_\g \,\chi_{_\rho}(\g)\, 2^{m-k+1}\, \frac{||\tf 12 f||^s}{(2\pi)^s}
 \left(  \sum_{j=0}^\infty \frac 1{(j+\frac 14)^s} - \sum_{j=0}^\infty
 \frac 1{(j+\frac  34)^s}\right)\\
 & = & (-1)^{r} \sigma_\g \, \chi_{_\rho}(\g)\, 2^{m-k+1} \, \frac{||\tf 12 f||^s}{(2\pi)^s}
 \left( \zeta(s, \tfrac  14) - \zeta(s, \tfrac 34)\right)
\end{eqnarray*}
where
 $\zeta(s,\al)$ is the  Riemann-Hurwitz function, for $\al \in (0,1]$.
Now $\zeta(s,\al)$ extends to an everywhere holomorphic function except for a simple
 pole at $s=1$, with residue 1 (see \cite{WW}, 13.13)
hence formula (\ref{etaz2kmanif}) implies that $\eta(s)$ is everywhere holomorphic.

Furthermore, since $\zeta(0,\al)=\frac 12 -\al$,  by
taking limit as $s\rightarrow 0$ in the above expression we get
(\ref{etainv}).
Also, differentiating (\ref{etaz2kmanif}) and using that $\zeta'(0,a)=\log
\G(a)-\frac12\log(2\pi)$ (see \cite{WW}) we obtain (\ref{p.etadif}).
\end{proof}

\begin{rem} Note that if $(\rho,V) = (1,\C)$, then  for $\Z_2^k$-manifolds with $k\le m$, one has that $\eta(0) \in 2\Z$. In particular the $\eta$-invariant of any
$\Z_2$-manifold is an even integer. Indeed, in the asymmetric case, $\eta=\pm 2^{m-1} \in \Z$. In dimension $n=3$, for $F\simeq\Z_2$, the proposition gives $\eta=\sigma_\g$.
 Take $\G=\langle \g, L_\Ld \rangle$ where $\g=\left[\begin{smallmatrix} -1 & & \\ & -1 & \\ & & 1 \end{smallmatrix}\right] L_{\f {e_3}2}$ and $\Ld$ the canonical lattice in $\R^3$. Then $\man$ has asymmetric spectrum only for the two spin structures $\vep_+=(1,1,-1;e_1e_2)$ and $\vep_-=(1,1,-1;-e_1e_2)$. Then $\eta_{(\G,\vep_+)}=1$ while $\eta_{(\G,\vep_-)}=-1$ as in \cite{Pf}.
\end{rem}

\section{Dirac isospectral manifolds}
In this section we give examples of  twisted Dirac isospectral flat manifolds that
are pairwise non-homeomorphic to each other. In Examples 4.3, 4.4 and 4.5
we compare (twisted) Dirac isospectrality with other types of isospectrality, such
as 
 Laplace isospectrality on functions and on $p$-forms and length
isospectrality with and without multiplicities (see Introduction). Two
manifolds are $[L]$-isospectral ($L$-isospectral) if they have the same
$[L]$-spectrum ($L$-spectrum). Obviously, $[L]$-isospectrality implies
$L$-isospectrality.

As a consequence we will obtain the following results:

\begin{teo} \label{comparison}
(i) There are families $\mathcal{F}$ of pairwise non-homeomorphic
Riemannian manifolds mutually twisted Dirac isospectral that are
neither Lapla\-ce isos\-pec\-tral on functions nor
$L$-isospectral. Furthermore, $\mathcal{F}$ can be chosen so that:

(a) Every $M \in \mathcal{F}$  has (no) harmonic spinors. (Ex.\@
\ref{z2manifolds} (i)).

(b) All $M$'s  in $\mathcal{F}$ have the same $p$-Betti numbers for $1\le
p \le n$ and they are $p$-isospectral to each other for any $p$ odd.
(Ex.\@ \ref{z2manifolds} (ii)).

(ii) There are pairs of non-homeomorphic spin manifolds that are
$\Delta_{s,\rho}$-isospectral but not $D_\rho$-isospectral.  (Ex.\@
\ref{example2} (ii)).

(iii) There are pairs of spin manifolds that are $\Delta_p$-isospectral
for $0\le p \le n$ and also $[L]$-isospectral which are
$D_\rho$-isospectral, or not, depending on the spin structure.
 (Ex.\@ \ref{example3} (i)).

 (iv) There are pairs of spin manifolds that are
$D_\rho$-isospectral and  $\Delta_p$-isospectral for
$0\le p \le n$ which are $L$-isospectral but not $[L]$-isospectral. (Ex.\@
\ref{example3} (ii)).
\end{teo}

\begin{teo} \label{exponential}
There exists a family, with cardinality depending exponentially on
$n$ (or $n^2$), of pairwise non-homeomorphic K\"ahler Riemannian
$n$-manifolds  that are twisted Dirac isospectral to each other
for many different spin structures. (Ex.\@ \ref{expfamily}, Rem.\@
4.7).
\end{teo}

To construct the examples, it will suffice to work with flat $n$-manifolds of
 diagonal type, having holonomy group $\Z_2^k$, $k=1,2$ and $n-1$.

 A Bieberbach group $\G$ is said to be of {\em diagonal type} if
 there exists an orthonormal $\Z$-basis $\{e_1,\dots,e_n\}$ of
the lattice $\Lambda$ such that for any element $BL_b\in\Gamma$,
$Be_i=\pm e_i$ for $1\le i\le n$ (see \cite{MR4}). Similarly,
$M_\G$ is said to be of diagonal type, if $\G$ is so. If $\Gamma$
is of diagonal type, after conjugation of $\Gamma$ by an isometry,
it may be assumed that $\Lambda$ is the canonical lattice and also
that  for any  $\gamma=BL_b \in \Gamma$, one has that $b \in \frac
1{2} \Lambda$ (\cite{MR4}, Lemma 1.4). \sk

Let $(M_\G, \vep)$ be a spin $\Z_2^k$-manifold where $\G=\langle
\g_1,\ldots,\g_k, L_\Ld \rangle$ and let $\ld_1,\ldots \ld_n$ be a
$\Z$-basis of $\Ld$. If $\g=BL_b\in \G$ we will fix  a
distinguished (though arbitrary) element in  $\mu^{-1}(B)$,
denoted by $u(B)$. Thus, $\vep(\g) = \sigma \, u(B)$,
 where $\sigma \in \{\pm 1\}$ depends on $\g$ and on the choice of $u(B)$.

The homomorphism $\vep:\G\rightarrow \spin$ is determined by its
action on the generators $\ld_1,\ldots,\ld_n,\g_1,\ldots,\g_k$ of
$\G$. Hence, if we set $\delta_i:= \vep(L_{\ld_i})$, we may
identify the spin structure $\vep$ with the $n+k$-tuple
\begin{equation}\label{spinstructure}
(\delta_1,\ldots,\delta_n, \sigma_1 u(B_1),\ldots,\sigma_k u(B_k))
\end{equation}
 where $\sigma_i$ is defined by the equation $\vep(\g_i)=\sigma_i u(B_i)$, for $1 \le i \le k$.

\begin{ejem}[{$\Z_2$-manifolds}] \label{z2manifolds}
 Here we will give some large twisted
Dirac isospectral sets of $\Z_2$-manifolds.
   Put $J:=\left[
\begin{smallmatrix} 0 & 1 \\ 1 & 0 \end{smallmatrix} \right]$. For each
$0\leq j,h < n$, define
\begin{equation}
 B_{j,h}:=\text{diag}(\underbrace{J,\dots,J}_j,
 \underbrace{-1,\dots,-1}_h,\underbrace{1,\dots,1}_l)
\end{equation}
 where $n=2j+h+l$, $j+h\not=0$ and $l\geq1$.  Then $\bjh \in \on$,
 $\bjh^2=\I$ and $\bjh \in \son$ if and only if $j+h$ is even.
Let $\Ld=\Z e_1\oplus \ldots \oplus \Z e_n$ be the canonical
lattice of
 $\R^n$ and for $j,h$ as before define the groups
\begin{equation}
 \G_{j,h}:=\langle \bjh L_{\frac{e_n}{2}}, \Ld \rangle.
\end{equation}
We have that $(\bjh +\I)\frac{e_n}2 = e_n \in \Ld \smallsetminus (\bjh
+\I)\Ld$. Hence, by Proposition~2.1  in \cite{DM}, the $\G_{j,h}$ are
Bieberbach groups. In this way, if we set $M_{j,h}= \G_{j,h} \backslash
\R^n$, we have a family
\begin{equation}                                        \label{family}
\mathcal{F}=\{ M_{j,h}=\G_{j,h} \backslash \R^n \,:\, 0 \le j \le
[\tfrac{n-1}2], 0 \le h < n-2j, \,j+h\not=0 \}
\end{equation}
of compact flat manifolds  with holonomy group $F\simeq \Z_2$.
The family $\mathcal{F}$ gives a system of representatives for
the diffeomorphism classes of $\Z_2$-manifolds of dimension $n$
(see Proposition 4.2 in \cite{MP}).


We will make use the following result on the existence of spin
structures (see Proposition 4.2 in \cite{MP}, and also Theorem 2.1
and Lemma 3.1):
\smallskip

{\em $M_{j,h}$ has $2^{n-j}$ spin structures parametrized by the
tuples $(\delta_1, \dots,\delta_n, \sigma ) \in \{\pm 1\}^{n+1}$
satisfying:}
 \begin{equation}                                    \label{relations}
   \delta_1=\delta_2, \;\cdots, \; \delta_{2j-1}=\delta_{2j} \qquad
   \text{ and } \qquad  \delta_n= (-1)^{\frac{j+h}2}.
 \end{equation}

\subsubsection*{Isospectrality on $p$-forms}
 We first review from \cite{MR3} and
 \cite{MR4} some results on  spectra of
 Laplace operators on vector bundles over flat manifolds. If $\tau$ is
a finite dimensional representation of $K=\on$ and $G=I(\R^n)$ we form the
vector bundle $E_\tau$ over $G/K \simeq \R^n$ associated to $\tau$ and
consider the corresponding vector bundle $\G\backslash E_\tau$ over
$\G\backslash \R^n =M_\G$. As usual we denote by $\chi_{_\tau}$ and
$d_\tau$ respectively, the character and the dimension of $\tau$. Let
$-\Delta_\tau$ be the connection Laplacian on this bundle.

We recall from  \cite{MR4}, Theorem 2.1,  that the multiplicity of the
eigenvalue $4\pi^2 \mu$ of $-\Delta_\tau$ is given by
\begin{equation}                                \label{eq.multip}
 d_{\tau,\mu}(\Gamma)= \tfrac{1}{|F|} \sum_{\gamma=BL_b
 \in \Lambda\backslash \Gamma}\chi_{_\tau}(B)
\; e_{\mu,\gamma} \;\;\text{ where }\;\; e_{\mu,\g} =
\sum_{{\scriptsize \begin{array}{c} v\in (\Ld^*)^B \\ \|v\|^2 =\mu
\end{array} } } e^{-2\pi i v\cdot b}.
\end{equation}

If $\tau=\tau_p$, the $p$-exterior representation of $\text{O}(n)$, then
$-\Delta_{\tau_p}$ corresponds to the Hodge Laplacian acting on $p$-forms.
In this case we shall write $\Delta_p$, $\text{tr}_p(B)$ and
$d_{p,\mu}(\G)$ in place of $\Delta_{\tau_p}$ , $\text{tr}\,\tau_p(B)$ and
 $ d_{\tau_p,\mu}(\G)$, respectively.

Thus, the $p$-Laplacian $\Delta_p$, $0\leq p\leq n$, has
eigenvalues $4\pi^2\mu$ with multiplicities $d_{p,\mu}$ given by formula
(\ref{eq.multip}). Furthermore, for flat manifolds of diagonal type the traces
$\text{tr}_p(B)$ are given by integral values of the {\em Krawtchouk
polynomials} $K_p^n(x)$ (see \cite{MR3}, Remark~3.6, and also \cite{MR4}).
Indeed, we have:
\begin{equation}                                \label{eq.krawtch}
    \text{tr}_p(B)=K_p^n(n-n_B), \quad
    \text{ where } K_p^n(x):= \sum_{t=0}^p (-1)^t \binom xt \binom {n-x}{p-t}.
\end{equation}

Note that $\bjh$ and $B_{0,j+h}$ are conjugate in $\text{GL}(n,\R)$ hence
$\text{tr}_p(\bjh)=\text{tr}_p(B_{0,j+h})=K_p^n(j+h)$. Thus,
\begin{equation}    \label{dpmu}
d_{p,\mu}(\G_{j,h})=\tfrac12 \big(\tbinom{n}{p} |\Ld_{\sqrt \mu}|
+ K_p^n(j+h) \,e_{\mu,\g}(\G_{j,h})\big ).
\end{equation}

Hence, the existence of integral zeros of
$K_p^n(x)$ will imply $p$-isospectrality of some of the $\mjh \in \fa$.
 For some facts on integral zeros of $K_p^n(x)$
see \cite{KL}, p.\@ 76, or Lemma 3.9 in \cite{MR4}. The simplest are the
so called trivial zeros, namely:
\sk

{\em If $n$ is even then $K_{\frac n2}^n(j)= K_j^n (\frac n2)=0$ for any
$j$ odd. Also,  $K_k^n(j)=0$ if and only if $K_j^n(k)=0$. } \sk

Thus for $n=2m$, all manifolds $\{\mjh: j+h \text{ is odd }\}$ are
$m$-isospectral and all  manifolds $\{\mjh: j+h=m \}$ are $p$-isospectral
for any $p$ odd with $1\leq p\leq n$. However, generically --i.e.\@, for
arbitrary $p$ and $n$--  the $\mjh$ will not be $p$-isospectral to each
other because the integral roots of Krawtchouk polynomials, aside from the
trivial zeros, are very sporadic.
\sk

We claim that the manifolds in $\fa$ are pairwise not isospectral on
functions. Take $\mu=1$. Then $\Ld_1=\{\pm e_1,\dots,\pm e_n\}$ and
$\Ld_1^{\bjh} = \{\pm e_{2j+h+1},\dots,\pm e_n\}$ thus $|\Ld_1|=2n$ and
$|\Ld_1^{\bjh}|=2(n-(2j+h))=2l$. Now, one checks that
$e_{1,\g}(\G_{j,h})=2(l-1)+2(-1)=2(l-2)$ and hence, from (\ref{dpmu}), we
get
\begin{equation} \label{dp1}
d_{p,1}(\G_{j,h})= \tbinom{n}{p} n + K^n_p(j+h)(l-2).
\end{equation}

Now consider $\mu=\sqrt 2$. Then
 $\Ld_{\sqrt 2}=\{\pm (e_i \pm e_j) : 1\le i < j \le n\}$ and
 $\Ld_{\sqrt 2}^{B_{j,h}}=\{\pm (e_{2i-1} + e_{2i}) : 1\le i \le j\}\cup
 \{ \pm(e_i\pm e_j) : l+1\le i < j \le n\}$. Hence
 $|\Ld_{\sqrt 2}|=4\tbinom n2$ and $|\Ld_{\sqrt 2}^{\bjh}|=2j+4 \tbinom l2$.
One checks that $e_{\sqrt 2,\g}(\G_{j,h})=2j - 4(l-1) +  (4\tbinom l2
-4(l-1))=2j+2(l-1)(l-4)$. Again by (\ref{dpmu}), we get
\begin{equation} \label{dp2}
d_{p,\sqrt 2}(\G_{j,h})= 2\tbinom{n}{p} \tbinom n2 +
K^n_p(j+h)\big(j+(l-1)(l-4)\big).
\end{equation}

In particular for $p=0$, since $K^n_0(j)=1$ for any $j$, we have
\begin{eqnarray}\label{d1}
d_{0,1}(\G_{j,h})& = & n + l - 2 \\
\label{d2} d_{0,\sqrt 2}(\G_{j,h})& =& n(n-1) + j + (l-1)(l-4).
\end{eqnarray}

These multiplicities are sufficient to show that all $\Z_2$-manifolds in
$\mathcal{F}$ are pairwise not isospectral. Indeed, if  $M_{j,h},
M_{j',h'}$ are isospectral then $l=l'$ by (\ref{d1}), thus $2j+h=2j'+h'$.
By (\ref{d2}), then $j=j'$ and hence $h=h'$.

\subsubsection*{Dirac isospectrality} 
We will give families of spin $\Z_2$-manifolds Dirac isospectral to each other.

We need to restrict ourselves to orientable manifolds, so
consider the family
$$\fa^+=\{ \mjh \in \fa : j+h \text{ is even}\}.$$
It will also be convenient to split $\fa^+=\fa^+_0 \dot \cup
\fa^+_1$ where 
$$\fa_i^+ = \{ \mjh \in \fa^+ : j+h\equiv 2i\mod(4) \}, \qquad i=0,1.$$

 (i) We now define spin structures for $\mjh$ in $\fa^+$. By
(\ref{relations}), $\delta_n=1$ for $j+h\equiv0 \, (4)$ and $\delta_n=-1$ for
$j+h\equiv2 \,(4)$. Hence, we take the spin structures
\begin{equation}                        \label{spinstructures}
\vep_{i,j,h}=(1,\ldots,1,(-1)^i;\sigma u(B_{j,h})), \qquad i=0,1,
\end{equation}
 for manifolds in $\fa^+_i$, $i=0,1$, 
 respectively. For simplicity,  we will write $\vep_{i}$ for
$\vep_{i,j,h}$.

We claim that all the spin $\Z_2$-manifolds in $\tilde{\fa}^+_0 := \{
(\mjh,\vep_{0}) : \mjh \in \fa^+_0 \}$ are twisted Dirac isospectral to
each other. Indeed, since $\vep_{0}$ is a spin structure of trivial type,
we know from Theorem~\ref{main2} that the spectrum is symmetric and the
multiplicities of the eigenvalues $\pm2\pi\mu$ of $D_\rho$ are given by
$$ d_{\rho,\mu}^{\pm}(\G_{j,h},\vep_0)=2^{m-2} d_\rho \,
|\Ld_{\vep_{0},\mu}^\ast| =2^{m-2} d_\rho\, |\Ld_{\mu}|$$ since
$\Ld^\ast_{\vep_{0}}=\Ld$. Note that all manifolds in $\tilde{\fa}^+_0$ have
$2^{m-1} d_\rho$ non-trivial harmonic spinors.

If $n\not \equiv 3 (4)$, the spin manifolds in $\tilde{\fa}^+_1 := \{ (\mjh,\vep_1) : \mjh \in
\fa^+_1 \}$  are Dirac isospectral to each other. The same happens  with those in $\tilde{\fa}^+_1 \smallsetminus \{ M_{0,2m} \}$, for $n=2m+1\equiv 3 (4)$.  Indeed, in both cases,  we have that
$d_{\rho,\mu}^{\pm}(\G_{j,h},\vep_1) = 2^{m-2} d_\rho\,
|\Ld_{\vep_1,\mu}|$, by Theorem~\ref{main2}. These manifolds do not have non-trivial harmonic spinors.
\sk

(ii) Note that, for every $t$, all $M_{t,0}, M_{t-1,1},\ldots, M_{0,t}$
have the same first Betti number. We recall from \cite{MP}, Proposition 4.1, that for $1\le
p\le n$
\begin{equation}       \label{bettip}
 \beta_p(M_{j,h})=  \sum _{i=0}^{[\frac p2]}\binom{j+h}{2i}
 \binom{j+l}{p-2i}.
\end{equation}
Hence, if $\beta_1(M_{j,h})=\beta_1(M'_{j,h})$, then
$\beta_p(M_{j,h})=\beta_p(M'_{j,h})$ for any $p \ge 1$.

Now, take
\begin{equation}
\mathcal{F}_t=\{(\mjh,\vep) : \mjh \in \fa^+ \text{ and } j+h=t \}
\end{equation}
for some fixed $t$ even  and $\vep$ as in (\ref{spinstructures}). In this
way $\mathcal{F}_t$ is a family of $t+1$ spin $\Z_2$-manifolds which are
Dirac isospectral to each other all having the same $p$-Betti numbers for
all $1\le p \le n$. Moreover, if we take $n=2t$ then
all $M_{j,h} \in \mathcal{F}_{t}$ are $p$-isospectral for any $p$ odd, by
the comments after (\ref{dpmu}).

\end{ejem}

\begin{ejem} \label{example2}
Here we give a simple pair of non-homeomorphic spin $\Z_2^2$-manifolds
that are (twisted) spinor Laplacian isospectral but not (twisted) Dirac isospectral. Let $\Ld$ be the canonical lattice in $\R^7$ and take the Bieberbach groups
$$\G=\langle B_1L_{b_1}, B_2L_{b_2}, \Ld \rangle, \qquad
\G'=\langle B_1L_{b_1'}, B_2L_{b_2'}, \Ld \rangle$$ where
$B_1=diag(-1,-1,-1,-1,-1,-1,1)$, $B_2=diag(-1,-1,1,1,1,1,1)$,
$B_1'=diag(-1,-1,-1,-1,1,1,1)$, $B_2'=diag(1,1,-1,-1,-1,-1,1)$, and
$b_1=\frac{e_7}2$, $b_2=\frac{e_1+e_3+e_7}2$, $b_1'=\frac{e_7}2$, $b_2'
=\frac{e_2}2$ are in $\tfrac12\Ld$. Let  $M_\G=\G\backslash\R^7$,
$M_{\G'}=\G'\backslash\R^7$ be the associated $\Z_2^2$-manifolds.

By Theorem 2.1 in \cite{MP}, one can check that $\man$ and $M_{\G'}$
respectively admit spin structures $\vep, \vep'$ with characters
$\delta_\vep=(\delta_1,\delta_2,\delta_1,\delta_4,\delta_5,\delta_6,-1)$ and
$\delta_{\vep'}=(\delta'_1,1,\delta'_3,\delta'_4,\delta'_5,\delta'_6,1)$,
$\delta_i, \delta_i' \in \{\pm1\}$.

Now, $F_1(\G)=\{B_1\}$. If we take $\delta_\vep=(1,1,1,1,1,1,-1)$, and if $\chi_\rho(B_1)\not=0$, then $(\man,\vep)$
has asymmetric Dirac spectrum. Thus, if $\mu_j=j+\frac 12$, $j\in \N_0$,
the multiplicity of $\pm 2\pi\mu$  is given by
$$ d_{\rho,\mu}^\pm(\G,\vep) =\left\{
\begin{array}{ll} d_\rho\,|\Ld_{\vep,\mu_j}|\pm 2(-1)^j \sigma_{\g_1} \chi_{_\rho}(\g_1) & \quad \mu=\mu_j
\sk \\ d_\rho\, |\Ld_{\vep,\mu}| & \quad \mu \ne \mu_j.
\end{array} \right.$$
 On the other hand $F_1(\G') = \emptyset$. Thus, $M_{\G'}$ has
symmetric Dirac spectrum with 
$d_{\rho,\mu}^\pm(\G',\vep')= d_\rho\, |\Ld_{\vep',\mu}|.$

Now take $\delta_{\vep'}=(-1,1,1,1,1,1,1)$. Then $(\man,\vep)$ and $(M_{\G'},\vep')$ are $\Delta_{s,\rho}$-isospectral. Indeed,
$d_0(\G,\vep)=d_0(\G',\vep')=0$ and since $|J_\vep^-|=|J_{\vep'}^-|=1$, for any $\mu>0$,  by Remark \ref{ortolattice}, we have
 $$d_{\rho,\mu}(\G,\vep)= 2 d_\rho\,|\Ld_{\vep,\mu}|=  2 d_\rho\, |\Ld_{\vep',\mu}|
 = d_{\rho,\mu}(\G',\vep').$$

However,  $(\man,\vep)$ and $(M_{\G'},\vep')$ are  not
$D_\rho$-isospectral for $Spec_{D_\rho}^A(\G,\vep)\not=\emptyset$ while
$Spec_{D_\rho}^A(\G',\vep')=\emptyset$. Note that $Spec_{D_\rho}^S(\G,\vep)=Spec_{D_\rho}^S(\G',\vep')$.
\end{ejem}

\begin{ejem}    \label{example3}
Here, we shall give two pairs of $\Delta_p$-isospectral for $0\le
p\le n$ and $L$-isospectral 4-dimensional $\Z_2^2$-manifolds
$M,M'$, one pair being $[L]$-isospectral and the other not. These
pairs will be twisted Dirac isospectral, or not, depending on the
choices of the spin structures.

Consider the manifolds $M_i,M_i'$, $1\le i \le 2$, where
$M_i=\G_i\backslash \R^4$, $M_i'=\G_i'\backslash \R^4$  and
$\G_i=\langle \g_1, \g_2, \Ld  \rangle$, $\G_i'=\langle \g_1',
\g_2', \Ld  \rangle$ are as given in Table 2, where $\g_i=B_i
L_{b_i}$, $\g_i'=B_i L_{b_i'},$  $i=1,2$, $B_3=B_1B_2$, $b_3=B_2
b_1 + b_2$, $b_3'=B_2' b_1' + b_2'$ and $\Ld=\Z e_1 \oplus \cdots
\oplus \Z e_4$ is the canonical lattice. Furthermore, we take $B_i
= B'_i$. In all cases the matrices $B_i$ are diagonal and are
written as column vectors. We indicate the translation vectors
$b_i, b'_i$ also as column vectors, leaving out the coordinates
that are equal to zero. We will also use  the pair
$\tilde{M_2},\tilde{M_2'}$ of $\Z_2^2$-manifolds of dimension 6
obtained from the pair $M_2,M_2'$  by adjoining the characters
$(-1,1,-1)$ and $(1,-1,-1)$ to $B_i$, $1 \le i \le 3$, and keeping
$b_i$, $b'_i$ unchanged.

Note that $M_1,M_1'$ and $\tilde{M}_2,\tilde{M}_2'$ are non-orientable while $M_2,M_2'$ not.

\renewcommand{\arraystretch}{1}
\medskip
\begin{center}
{\sc Table 2}

\smallskip
 $ \left. \begin{array}{c} \{M_1,M_1'\} \end{array} \right.$ \qquad
\begin{tabular}{|rcc|rcc|rcc|}  \hline  $B_1$ &  $L_{b_1}$  & $L_{b_1'}$    & $B_2$ &  $L_{b_2}$ & $L_{b_2'}$  & $B_3$ &  $L_{b_3}$ & $L_{b_3'}$   \\ \hline
-1 & & &      1 & &  $\text{{\scriptsize 1/2}}$ & -1 & &$\text{{\scriptsize 1/2}}$ \\
-1 & & &   -1 & $\text{{\scriptsize 1/2}}$ & $\text{{\scriptsize 1/2}}$& 1 &$\text{{\scriptsize 1/2}}$ & $\text{{\scriptsize 1/2}}$    \\
1  & & $\text{{\scriptsize 1/2}}$  &      -1 && & -1 & &  $\text{{\scriptsize 1/2}}$   \\
1  & $\text{{\scriptsize 1/2}}$ &  & 1  & $\text{{\scriptsize 1/2}}$ & & 1  & & \\
\hline
\end{tabular}

\smallskip
$\left.\begin{array}{c} \{M_2,M_2'\}\\
\{\tilde{M}_2,\tilde{M}'_2\}
\end{array}\right.$
 \qquad
\begin{tabular}{|rcc|rcc|rcc|}  \hline    $B_1$ &  $L_{b_1}$  & $L_{b_1'}$    & $B_2$ &  $L_{b_2}$ & $L_{b_2'}$  & $B_3$ &  $L_{b_3}$ & $L_{b_3'}$   \\ \hline
 1 & & & 1 & & $\text{{\scriptsize 1/2}}$ & 1 & & $\text{{\scriptsize 1/2}}$ \\
 1 & & $\text{{\scriptsize 1/2}}$ & 1 & $\text{{\scriptsize 1/2}}$ & $\text{{\scriptsize 1/2}}$ & 1 &$\text{{\scriptsize 1/2}}$ &  \\
 1 &  $\text{{\scriptsize 1/2}}$ & &  -1 & & & -1 &  $\text{{\scriptsize 1/2}}$ &  \\
 -1 &  & & 1  &$\text{{\scriptsize 1/2}}$ & & -1  & $\text{{\scriptsize 1/2}}$  & \\
 \hline
-1 &  & &  1 & & & -1 &  &  \\ 1 &  & & -1 & & & -1 & & \\  \hline
\end{tabular}
\end{center}



\medskip
(i) By Example 3.3 in \cite{MR4}, the manifolds $M_1, M_1'$ are Sunada
isospectral (hence $p$-isospectral for $0\le p\le n$) and
$[L]$-isospectral.

By Theorem 2.1 in \cite{MP}, one can check that $M_1$ admits $2^4$ spin
structures $\vep_1$ of the form
$\vep_1=(\delta_1, -1, \delta_3, -1; \sigma_1 e_1 e_2, \sigma_2 e_2e_3)$,
and  $M_1'$ carries $2^3$ spin structures $\vep_1'$ of the form
$\vep_1'=(-1, -1, -1, \delta_4', \sigma_1' e_1 e_2, \sigma_2' e_2 e_3)$
where $\delta_1,\delta_3,\delta_4',\sigma_1,\sigma_2,\sigma_1',\sigma_2'
\in \{\pm 1\}$.

Since $F_1(\G_1)=F_1 (\G'_1)=    \emptyset$, by (\ref{multeven})
we have that  $d_{\rho,\mu}^\pm(\G_1,\vep_1)= d_\rho
\,|\Ld_{\vep_1,\mu}|$ and $d_{\rho,\mu}^\pm(\G_1',\vep_1')=
d_\rho\, |\Ld_{\vep_1',\mu}|$.

Now, if we take $\delta_{\vep_1}=(1,-1,1,-1)$, we see that
$(M_1,\vep_1)$ is not twisted Dirac isospectral to
$(M_1',\vep_1')$ for any $\vep_1'$ because $|J_{\vep_1}^-|=2$
while $|J_{\vep_1'}^-|\geq3$ (see Remark \ref{ortolattice}).
However, if we take $\vep_1, \vep_1'$ such that
$\delta_{\vep_1}=\delta_{\vep_1'}=(-1,-1,-1,-1)$, then
$(M_1,\vep_1)$ and $(M_1',\vep_1')$ \! are $D_\rho$-isospectral to
each other\!.

\medskip
(ii) By Example 3.4 in \cite{MR4}, the manifolds $M_2, M_2'$ are
Sunada isospectral (hence $p$-isospectral for $0\le p \le n$ and
$L$-isospectral) but not $[L]$-isospectral. In order to have
orientable manifolds we add to $M_2,M_2'$ the characters
$(-1,1,-1)$ and $(1,-1,-1)$. The pair $\tilde{M}_2,\tilde{M}_2'$
obtained has the same spectral properties as $M_2,M_2'$. This can
be seen by proceeding as in  Example 3.4 in \cite{MR4}.

Again, by Theorem 2.1 in \cite{MP}, we can check that $\tilde{M}_2$ has
$2^5$ spin structures, $\vep_2$, with characters
$\delta_{\vep_2}=(\delta_1,1,-1,-1,\delta_5,\delta_6)$,
$\delta_1,\delta_5,\delta_6 \in\{\pm1\}$; and $\tilde{M}_2'$ has $2^6$
spin structures, $\vep_2'$, with characters
$\delta_{\vep_2'}=(1,-1,\delta_3',\delta_4',\delta_5',\delta_6')$ where
$\delta_3',\delta_4',\delta_5',\delta_6' \in \{\pm1\}$.

Now,  $d_{\rho,\mu}^\pm(\tilde{\G}_2,\vep_2)=
d_\rho\,|\Ld_{\vep_2,\mu}|$ and
$d_{\rho,\mu}^\pm(\tilde{\G}_2',\vep_2')= d_\rho\,
|\Ld_{\vep_2',\mu}|$. As in (i), if we take spin structures
$\vep_2, \vep_2'$ with $|J_{\vep_2}^-|=|J_{\vep_2'}^-|$, we see
that $(\tilde{M}_2,\vep_2)$ and $(\tilde{M}_2',\vep_2')$ are
$D_\rho$-isospectral to each other, whereas if we take $\vep_2,
\vep_2'$ such that
$|J_{\vep_2}^-|\not =|J_{\vep_2'}^-|$, then $(\tilde{M}_2,\vep_2)$
and $(\tilde{M}_2',\vep_2')$ 
are not $D_\rho$-isospectral.

\end{ejem}

\begin{ejem}    \label{expfamily}
We now construct a large family of pairwise non-homeomor\-phic
twisted Dirac isospectral flat $2n$-manifolds, with holonomy group $\Z_2^{n-1}$.
We will apply the doubling procedure in \cite{JR} or in \cite{BDM} to the family of Hantzsche-Wendt manifolds (see \cite{MR}).

  We first recall some facts from \cite{MR}.
Let $n$ be odd. A {\it Hantzsche-Wendt  group} (or {\it HW   group})
 is an $n$-dimensional orientable Bieberbach group $\Gamma$ with holonomy
 group $F\,\simeq\,\Z_2^{n-1}$ such that  the action of every $B\in F$
 diagonalizes on the canonical  $\Z$-basis $e_1,\dots, e_n$  of
 $\Lambda$. The holonomy group $F$
 can thus be identified to  the diagonal subgroup
$\{B:Be_i = \pm e_i, \,\,1\,\le\,i\,\le\,n,\,\,\det B = 1\}$ and
$M_\Gamma =\Gamma\backslash\R^n$ is called a {\it Hantzsche-Wendt}
 (or {\it HW}) {\it  manifold}.

We denote by $B_i$ the diagonal matrix  fixing $e_i$ and such that
$B_i e_j = -e_j$ (if $j \ne i)$, for each $1\,\le\,i\,\le\,n$.
Clearly, $F$ is generated by $B_1, B_2,\dots, B_{n-1}$ and furthermore
$B_n = B_1 B_2 \dots B_{n-1}$.

Any HW group has the form $\Gamma = \langle \g_1, \dots,\g_{n-1},
L_\lambda\,:\,\lambda\,\in\, \Lambda\rangle$, where $\g_i = B_i
L_{b_i}$ for some $b_i\in\R^n$, $1\le i \le n-1$, and where it may
be assumed that the components $b_{ij}$ of $b_i$ satisfy
$b_{ij}\in\{0,\frac12\}$, for $1\,\le\,i,j\,\le\,n$. 
Also, since $({\Ld^p(\R^n)})^{\text{F}}=0$ for any $1\le p \le
n-1$, it follows that all HW manifolds are rational homology
spheres. We further recall that, as shown in \cite{MR}  (by using
a rather small subfamily ${\mathcal H}_1$), the cardinality $h_n$
of the family of all  HW groups under isomorphism satisfies $h_n
>\frac {2^{n-3}}{n-1}$. Moreover, the cardinality of the set of
Laplace isospectral, pairwise non-isomorphic, {\it pairs} of HW
groups grows exponentially with $n$.
\sk

Relative to the spin condition, it is easy to verify using
($\varepsilon_1$) in \cite{MP}, (2.3), that the manifolds in this
family are generally not spin. Indeed, we now show that no HW
manifold with $n=4k+1$ is spin. We note that $\g_k^2= L_{e_k}$ for
each $1\le k\le n$. Hence a spin structure $\vep$ must satisfy
$\delta_k =\vep(L_{e_k}) = \vep(\g_k^2) =\vep(\g_k)^2= (\pm
e_1\ldots e_{k-1} e_{k+1} \ldots e_n)^2 = 1$, since  $n = 4k+1$.
Thus,  it follows that $\vep_{|\Ld}=\I$. On the other hand $(\g_i
\g_j)^2 = L_\ld$ for some $\ld \ne 0$ and since $\vep (\g_i \g_j)
=\pm e_i e_j$ it follows that $\vep(L_\ld) = \vep  (\pm e_i e_j)^2
=-1$, a contradiction. \sk

 In the case $n= 4k+3$, $k>0$, we shall
show that no HW manifold in the family ${\mathcal H}_1$ (see
\cite{MR}) is spin. More generally, assume that there are three
consecutive generators of $\G$, $\g_i = B_iL_{b_i}$ with $b_i
=L_{(e_i +e_{i-1})/2}$, $i\ge 2$. Thus, if
$\g:=\g_i\g_{i+1}\g_{i+2}= B_i B_{i+1}B_{i+2}L_{(e_{i-1}
+e_{i+2})/2}$, then we have $\g^2 = \g_{i+2}^2=L_{e_{i+2}}$. This
gives a contradiction since $\vep(\g)^2=1$ (the multiplicity of
the eigenvalue $-1$ for $B_iB_{i+1}B_{i+2}$ is $4k$) while
$\vep(\g_{i+2})^2 =-1$ (the multiplicity of the eigenvalue $-1$
for $B_{i+2}$ is $4k+2$).
Actually, J.P. Rossetti has recently shown to us a proof (that still uses the criterion in \cite{MP}) that no HW manifold can admit a spin structure.

\sk

We will now consider for any  HW group $\G$, the group $\dd\G$, defined by the doubling construction in \cite{JR} or \cite{BDM}, namely $\dd\G =\langle dB L_{(b,b)},\,
L_{(\ld_1,\ld_2)} : BL_b\in \G, \ld_1, \ld_2 \in \Ld \rangle$ where $d B :=\left[ \begin{smallmatrix} B & 0 \\
0 & B  \end{smallmatrix} \right]$.  This yields a Bieberbach group of
dimension  $2n$, with the same holonomy group $\Z_2^{n-1}$ as $\G$, and
with the additional property that the associated manifold $M_{\dd\G}$ is
K\"ahler. The reason why we use ${\dd\G}$ in place of $\G$ is that  $M_{\dd\G}$ is always spin  (see \cite{MP}, Remark 2.4).

We will need the following facts:

\msk
{\em (i) If $\G$ is a HW group, then  $M_{\dd\G}$ admits $2^{n-1}$
  spin structures of trivial type.}

{\em (ii) If $\G$ runs through all $HW$ groups, all  manifolds $M_{\dd\G}$ endowed with   spin structures of trivial type are twisted Dirac isospectral to each other.}

{\em  (iii) If $\G, \G'$ are non-isomorphic HW groups then $\dd\G, \dd\G'$
are non-\-iso\-mor\-phic Bieberbach groups.}

{\em  (iv) Two HW groups $\G, \G'$ are Laplace isospectral if and only if
$\dd\G$ and $\dd\G'$ are Laplace isospectral.}

\msk
Now, (i) and (ii) are direct consequences of Remark 2.4 in \cite{MP} and
of (i) of Theorem \ref{main2}, respectively.

\sk {\em Proof of (iii)}. This assertion follows by an argument very
similar to that given in the proof of Proposition 1.5 in \cite{MR}. We
shall only sketch it. An isomorphism between $\dd\G$ and $\dd\G'$ must be
given by conjugation by $CL_c$ with $C \in \text{GL}(2n, \R)$ and $c
\in\R^{2n}$. Now $C(\Ld \oplus \Ld) =\Ld \oplus \Ld$ implies $C\in
\text{GL}(2n, \Z)$ and furthermore for each $1\le i \le n-1$,
$$CL_cdB_iL_{(b_i,b_i)}L_{-c}C^{-1}= dB_{\sigma(i)}
L_{(b'(\sigma(i),b'(\sigma(i))}$$ where $\sigma \in S_n$. Thus
$CdB_iC^{-1}L_{C((b_i,b_i) + (dB_i -\I)c)}.$ In particular, this
implies that $CdB_iC^{-1}= dB_{\sigma (i)}$ for each $1\le i \le
n$, with $\sigma \in S_n$. Thus, there is an $n\times n$
permutation matrix $P$ such that $D:=C dPC^{-1}$ commutes with
$dB_i$ for each $i$, thus $D$ preserves $\Z e_i \oplus \Z
e_{n+i}$ for each $i$. It is easy to see that conjugation by
such $D$ yields an automorphism of $\dd\G'$. Thus, conjugation
by $DCL_c = dPL_c$ takes  $\dd\G$ onto $\dd\G'$ isomorphically
and furthermore
\begin{eqnarray*}dPL_c dB_iL_{(b_i,b_i)}L_{-c}dP^{-1}&=&
 d(PB_iP)^{-1}L_{dP((b_i,b_i) + (dB_i -\I)c)}
 \\& = & dB_{\sigma (i)}L_{(b'_i,b'_i)}.
\end{eqnarray*}
This implies that $c= (c_1,c_1)$, mod $\Ld$, with $c_1 \in \frac 14 \Ld$
and hence, conjugation by $PL_{c_1}$ gives an isomorphism between $\G$
and $\G'$,  a contradiction.

\sk
{\em Proof of (iv)}.
 Since HW groups are of diagonal type, then $\G$, $\G'$ are
Laplace isospectral if and only if they are Sunada isospectral, that is,
if they have the same Sunada numbers (see \cite{MR3}, \cite{MR4}). We
claim that this is the case if and only if $\dd\G$ and $\dd\G'$ are Sunada
isospectral to each other.

Indeed, we recall that for  $0\le t\le s \le n,$ and $\G$ of diagonal
type,  the {\it Sunada numbers} of $\G$ are defined by
\begin{equation*}\label{sunada}
c_{d,t}(\G) :=\big|\big\{BL_b \in \G : n_B=d \text{ and } n_B(\tfrac 12)
=t \big \}\big|.
\end{equation*}
where, for $BL_b \in \G$, $n_B:= \text{dim}(\R^n)^B=|\{ 1\le i \le
n: Be_i=e_i \}|$ and $n_B(\tfrac 12) := |\{1 \le i \le n: Be_i=e_i
 \text { and } b\cdot e_i \equiv \tfrac 12 \mod(\Z) \}|$.
 Now, it is clear from the definitions that
$c_{2s,2t}(\dd\G)= c_{s,t}(\G)$ for each $0\le t\le s \le n$ and
$c_{u,v}(\dd\G)= 0$, if either $u$ or $v$ is odd. This clearly implies
that $\G$ and $\G'$ have the same Sunada numbers if and only if $\dd\G$
and $\dd\G'$ do so, that is, if and only if $\dd\G$ and $\dd\G$ are
Sunada-isospectral to each other.

\msk
Thus, for each $n$ odd, by (i), (ii), (iii), (iv), the above construction
yields a family, of cardinality that depends exponentially on $n$, of K\"ahler flat manifolds of dimension $2n$,
all pairwise non-homeo\-mor\-phic and all twisted Dirac isospectral to each other,
 having only $2 d_\rho$ harmonic spinors for every trivial spin structure chosen (see (iii) of Theorem \ref{main2}).
Within this family, by (iv) and \cite{MR}, there are exponentially
many pairs that are Sunada isospectral, hence $p$-isospectral for
all $p$. However,  generically, two such  manifolds  will not be
$p$-isospectral for any value of $p$ (see for instance \cite{MR2}
in the case $n=7$).

We note that if we repeat the doubling procedure then the set of all
$M_{{\dd}^2\G}$, with $\G$ a HW group, is an exponential family of
hyperk\"ahler manifolds with the same spectral properties as the family
$M_{\dd\G}$, but now having $2^{n+1} d_\rho$ harmonic spinors for every
spin structure of trivial type chosen.
\end{ejem}

\begin{rem}
If one looks at the family of all flat manifolds with holonomy
group $\Z_2^{n-1}$ (see \cite{RS}) for $n=4r+3$, then it is shown
in \cite{MPR} that a subfamily of this family has cardinality
$2^{\frac{(n-1)(n-2)}2}$. If we apply the doubling procedure to
this family, one shows that considerations
$\text{(i)},\ldots,\text{(iv)}$ in Example 4.6 remain valid, hence
one obtains a family of twisted Dirac isospectral,
  pairwise non homeomorphic, $2n$-manifolds of cardinality  $2^{\frac{(n-1)(n-2)}2}$.
\end{rem}

\section{Eta invariants of $\Z_p$-manifolds}
In this section we shall illustrate the results in Section 2  by
using   the expression (\ref{etaodd}) of the eta series, to
compute explicitly the eta invariant of certain flat $p$-manifolds
with cyclic holonomy group $\Z_p$, $p=4r+3$, $p$  prime, for the
two existing spin structures. In \cite{SS}, the authors give an
expression for the eta invariant and harmonic spinors of this
family (without assuming $p$ to be prime), in terms of the
solutions of certain equations in congruences. They give explicit
values for $p=3,7$. Here we shall give an explicit expression for
the eta invariant of this family in terms of Legendre symbols and
special values of trigonometric functions. At the end we give a
table with the values of $\eta$ for any prime $p\le 503$. For
$n=3$ they coincide with those computed in \cite{Pf} (and
\cite{SS}). Our formulas for the eta series involve trigonometric
sums and resemble those obtained in \cite{HZ} to compute the
$G$-index of elliptic operators for certain low dimensional
manifolds.

We thus assume that $F$ is cyclic of order $p=4r+3$ (i.e.\@ $m=2r+1$) with
$p$ prime.
Let  $\G =\langle BL_b , L_\Ld \rangle$. Here $\Ld =\sum_{i=1}^p \Z e_i$,   $B$ is of order $p$ with $B(e_i)=e_{i+1}$ for $1\le i\le p-2$, $B(e_{p-1}) =-\sum_{i=1}^{p-1} e_i$, $B(e_p) = e_p$ and $b=\tf 1p e_p$.

Now, there exists $D \in \text{GL}_p(\R)$, $De_p =e_p$ such that $C:=DBD^{-1} \in \text{SO}(p)$. Thus, $\G_p:= D\G D^{-1}=\langle \g=CL_{\f {e_p}p}, D\Ld \rangle$ is a Bieberbach group and $M_p= \G_p \backslash \R^p$ is an orientable $\Z_p$-manifold of dimension $p$. The vectors $f_i :=De_i$ for $1\le i \le
p-1$, and $f_p=e_p$ give a $\Z$-basis of $D\Ld$.

Since the eigenvalues of $B$ and $C$ are the $p$-roots of unity:
$\{ e^{\f {2 \pi i k}{p}}: 1\le k \le p-1 \}$, then we can assume
by further conjugation in $\text{SO}(p)$ that $C=
x_0(\tfrac{2\pi}p,\ldots, \tfrac {2m\pi}p)$ (see
(\ref{torielements})) with $m= \tf {p-1}2$.

We now note that Theorem 2.1 in \cite{MP}, stated for
$\Z_2^k$-manifolds, also holds for $\Z_n$-manifolds, replacing
condition $(\vep_1)$: $\vep(\g^2)=\vep(\g)^2$,  by $(\vep_1')$:
$\vep(\g^n)=u_B^n$ for any $u_B \in \spin$ such that $\mu(u_B)= B$
and keeping condition $(\vep_2)$  in
\cite{MP}, i.e. $\vep(L_{(B-\I)(\ld)})=1$ for any $\ld \in \Ld$.
Thus, using conditions $(\vep_1')$ and $(\vep_2)$ one can see that
$M_p:= \Gamma_p\backslash \R^p$ has exactly two spin structures,
$\vep_1, \vep_2$, given, on the generators of $\G$, by
\vspace{-.150cm}
\begin{gather*}
\vep_1(L_{f_j})=\vep_2(L_{f_j}) =1 \,\,(1\le j \le p-1),\quad
  \vep_1(L_{e_p})= 1,\, \vep_2(L_{e_p}) =-1, \\
\vep_1 (\g)= (-1)^{r+1} x(\tf{\pi}p, \tf{2\pi}p,\dots, \tf{(2r+1)\pi}p),
\quad \vep_2 (\g)= (-1)^r x(\tf{\pi}p, \tf{2\pi}p,\dots, \tf{(2r+1)\pi}p),
\end{gather*}
in the notation of (\ref{torielements}). We can now state:

\begin{teo} Let $p=4r+3$ be a prime and let $\G_p$, $\vep_1, \vep_2$ be
as above. Then, the eta series of
$M_p=\G_p\backslash \R^p$, for $\vep_1$ and $\vep_2$ are respectively given by:
\begin{equation}\begin{split} \label{etas}
\eta_{(\G_p,\rho,\vep_1)}(s) & = \tfrac{-2\, \chi_\rho(\g)}{\sqrt p\,(2\pi p)^s} \,
 \sum_{k=1}^{p-1}  \, \Big( \f {k}p \Big)
\sum_{l=1}^{p-1}\sin(\tf{2 l\pi  k}p) \, \zeta(s, \tf lp), \\
\eta_{(\G_p,\rho,\vep_2)}(s) & = \tfrac{-2 \,\chi_\rho(\g)}{\sqrt p\,(2\pi p)^s} \;
\sum_{k=1}^{p-1}  \, (-1)^k \Big( \f {k}p \Big)
\sum_{l=0}^{p-1}\sin(\tf{(2l+1) \pi k}p) \, \zeta(s, \tf {2l+1}{2p}),
\end{split}
\end{equation}
and the eta invariants have the expressions
\begin{equation}  \begin{split}
\label{etainvariants}
\eta_{\rho,\vep_1}  & = \tf{-2 \,\chi_\rho(\g)}{\sqrt p} \; \sum_{k=1}^{\f{p-1}2}\,
\Big( \f kp \Big)\, \text{cot} (\tfrac {k\pi}p), \quad
\\
\eta_{\rho,\vep_2}  & = \tf{-2 \,\chi_\rho(\g)}{\sqrt p} \; \sum_{k=1}^{\f{p-1}2}\,
(-1)^k \Big( \f kp \Big) \, \text{cosec} (\tfrac {k\pi}p),
 \end{split}
\end{equation}
where  $\Big(\f \cdot p \Big)$ denotes the Legendre symbol.
\end{teo}
\noindent{\em Note.} The theorem shows that one can have a spin
structure of trivial type ($\vep_1$) with asymmetric Dirac
spectrum, in contrast with the situation in the case of holonomy
group $\Z_2^k$.

Expressions (\ref{etas}) and (\ref{etainvariants}) can be
simplified further by using identities in number theory. In
particular one shows that the eta invariants take integer values.
We plan to get deeper into this question in a sequel to this
paper.

\begin{proof} We will compute the
different ingredients in formula (\ref{etaodd}) for the eta
function for the two given spin structures.  For $1\le k\le p-1$
we have $(\Ld_{\vep}^\ast)^{B^k} = \Z e_p$ for $\vep_1$ while
$(\Ld_{\vep}^\ast)^{B^k} = (\Z + \frac 12) e_p$ for $\vep_2$. In
both cases $(\Ld_{\vep,\mu}^\ast)^{B^k} = \{\pm \mu e_p \}$. Thus,
$\mu=j$ for $\vep_1$ and $\mu=j - \tf 12$ for $\vep_2$ with $j\in
\N$. We take $x_{\g^k} = \vep_i(\g)^k$, for $i=1,2$. This implies
that $\sigma(e_p, x_{\g^k})=1$ for $\vep_1,\vep_2$.

 Since $\g^k =B^k L_{\f{k}{p}e_p}$,  by (\ref{emugamas}) and
Remark
\ref{sigmas}, for $1\le k \le p-1$, we have 
$$e_{\mu,\g^k,\sigma}(\delta_{\vep})= e^{-2 \pi i \f{\mu k}{p}} \sigma(\mu
e_p, x_{\g^k})  + e^{2 \pi i \f{\mu k}p} \sigma(- \mu e_p, x_{\g^k}) =-2i
\sin(\tf{ 2 \pi \mu k}p).$$
Thus, up to the factor $-2i$, the sums over $\mu \in \frac{1}{2\pi}
\mathcal{A}$ corresponding to  $\g=\g^k$ in (\ref{etaodd}) for $\vep_1$
and  $\vep_2$ are respectively equal to 
 \begin{gather}\label{sinszetas}  \nonumber
 \qquad   \quad \sum_{j=1}^\infty
  \frac{\sin(\f {2 j \pi  k}p)}{j^s}  =
   \sum_{l=1}^{p-1}\sin(\tf{2 l \pi k}p) \sum_{j=1}^\infty
   \frac{1}{(pj+l)^s}  = \f 1{p^{s}} \sum_{l=1}^{p-1}
   \sin(\tf{2 l \pi  k}p) \zeta(s,\tf lp),
 \\ \vspace{-.5cm}   \\ \nonumber
  \sum_{j=0}^\infty
  \frac{\sin(\f{(2j+1)\pi k}p)}{(j+\f 12)^s}  =
  \f{1}{p^{s}} \sum_{l=0}^{p-1} \sin(\tf{(2l+1) \pi k}p)
 \zeta(s,\tf{2l+1}{2p}),
\end{gather}
where $\zeta(s,\al) = \sum_{j=0}^\infty \frac 1{(j+\alpha)^s}$ denotes the
Riemann-Hurwitz zeta function for $\alpha \in (0,1]$.

\smallskip
We now compute the product of sines in (\ref{etaodd}) in both
cases. We note that  $t_j(x_{\g^k})$ (see (\ref{torielements}))
depends on $\g^k$ and also on $\vep_1, \vep_2$.
We have, for $\vep_h$, $h=1,2$: 
\begin{eqnarray*}
\Pi_{j=1}^{\f{p-1}2} \sin t_j(x_{\g^k})&=& (-1)^{k(r+h)}\,
\Pi_{j=1}^{\f{p-1}2} \sin (\tfrac{\pi j k}p)  \\
&=& (-1)^{k(r+h)}\, \Pi_{j=1}^{\f{p-1}2} (-1)^{[\f{jk}p]}\,
\Pi_{j=1}^{\f{p-1}
2} \sin (\tfrac{\pi j}p)\\&=&
 (-1)^{k(r+h)} (-1)^{s_p(k)} \frac{\sqrt p}{2^{2r+1}},
\end{eqnarray*}
where we have put $s_p(k)= \sum_{j=1}^{\f{p-1}2} [\frac{jk}p]$ and used
the identities
$$\sin(\pi z) =\frac \pi{\G(z)\G(1-z)}, \quad
(2\pi)^{\tfrac {p-1}2} \G(z)= p^{z-\frac 12} \G(\tfrac zp) \G(\tf
{z+1}p)\cdots   \G(\tf {z+p-1}p).$$
Now, if $\Big( \f {\cdot}p \Big)$ denotes the Legendre symbol, then  (see \cite{Ap}, Theorems
9.6 and 9.7) we have, since $p=4r+3$,
$$ (-1)^{s_p(k)}=(-1)^{(k-1)(\f {p^2-1}8)} \Big( \f {k}p \Big) =
(-1)^{(k-1)(r+1)} \Big( \f {k}p \Big).$$
In this way, we obtain:
\begin{equation} \label{productsines}
\Pi_{j=1}^{\f{p-1}2} \sin t_j(x_{\g^k}) = \left\{ \begin{array}{ll}
(-1)^{(r+1)}\, 2^{-2r-1} \sqrt{p} \Big( \f {k}p \Big)  & \quad \text{
for } \vep_1  \\ (-1)^{(r+1)}
 \, 2^{-2r-1} (-1)^k \sqrt{p}\Big( \f {k}p \Big) 
& \quad \text{ for } \vep_2.
\end{array} \right.
\end{equation}

Now, starting from (\ref{etaodd}) and using (5.3) 
and (\ref{productsines}) we finally arrive at the expressions for the
$\eta$-series of $(M_p, \vep_h)$, $h=1,2$, given in (\ref{etas}).
\sk

We now compute the eta invariants.
Using that $\zeta(0,\al)= \f 12-\al$ (\cite{WW}, 13.21), together with
the fact that  $\sum_{l=1}^{p-1}\sin(\tf{2l \pi  k}p)=
\sum_{l=0}^{p-1}\sin(\tf{(2l+1) \pi k}p)=0$, for every $1\le k \le p-1$,
we see that
the sums over $l$ in the expressions (\ref{etas}),
when evaluated at $s=0$ are respectively equal to
\begin{equation}\label{sums}
-\f 1{p}\sum_{l=1}^{p-1}
 l\sin(\tf{2l \pi  k}p), \text{ for } \vep_1, \quad
-\f 1{2p}\sum_{l=0}^{p-1} (2l+1)\sin(\tf{(2l+1) \pi k}p),
\text{ for } \vep_2.
\end{equation}
We claim that
\begin{equation}\label{sinesums}
\sum_{l=1}^{p-1}\,l \sin(\tf{2l \pi  k}p)= -\tfrac p2 \text{cot} (\tfrac
{k\pi}p),\quad
 \sum_{l=0}^{p-1}\, l\sin\big(\tf{(2l+1)\pi k}p\big)= -\tfrac p2 \text{cosec}(\tfrac {k\pi}p).
\end{equation}
Indeed, by differentiating the identity $\tf 12 + \sum _{l=1}^{p-1}
\cos(lx) = \frac{\sin ((p-\f 12)x)}{2\sin (\tf x2)}$, we get
$$
-\sum_{l=1}^{p-1}\, l \sin (lx) =\frac {(2p-1) \cos((p-\tf 12)x)}{4\sin
(\tf x2)} - \frac {\cos (\tf x2) \sin((p-\tf 12)x)}{4\sin^2 (\tf x2)}.
$$
Evaluating both sides at $x=\tf{2k\pi}p$ yields the first equality in
(\ref{sinesums}).

To verify the second identity  we first note that
$$
\sum_{l=1}^{p-1} \, l \cos (\tf {2l\pi k}p)= p\sum_{l=1}^{\f{p-1}2} \,
\cos (\tf {2l
\pi k}p)= -\tf p2.
$$
 Using this expression together with the first identity in (\ref{sinesums}) we have
\begin{eqnarray*}
 \sum_{l=0}^{p-1}\, l\sin\big(\tf{(2l+1)\pi k}p\big)&=&
 \cos (\tf{\pi k}p) \text{cot} (\tf{\pi k}p)(-\tf p2 )
 +\sin (\tf{\pi k}p) (-\tf p2) 
\\&=&
- \tf p2 \text{cosec}(\tf{\pi k}p).
\end{eqnarray*}

Hence, from (\ref{etas}), (\ref{sums}) and  (\ref{sinesums}) we obtain
$$ \eta_{\vep_1} = \f{-1}{\sqrt p} \, \sum_{k=1}^{p-1}\, \Big( \f
{k}{p} \Big)\, \text{cot} (\tfrac {k\pi}p), \quad
 \eta_{\vep_2} =
\f{-1}{\sqrt p} \, \sum_{k=1}^{p-1}\, (-1)^k \Big( \f {k}{p} \Big)\,
\text{cosec} (\tfrac {k\pi}p).$$ We finally note that the contributions of $k$
and $p-k$ to the above sums are equal to each other, that is:
\begin{gather*}
  \Big( \f {k}{p} \Big) \text{cot} (\tfrac {k\pi}p) =
   \Big( \f {p-k}{p} \Big)   \text{cot} (\tfrac {(p-k)\pi}p)  \\
(-1)^k \Big( \f {k}{p} \Big)  \text{cosec} (\tfrac {k\pi}p) = (-1)^{p-k}
\Big(\f {p-k}{p} \Big) \text{cosec} (\tfrac {(p-k)\pi}p).
\end{gather*}
These identities can be easily verified using that
$$ \Big( \f {p-k}{p} \Big) = \Big( \f {-k}{p} \Big) =
\Big( \f {-1}{p} \Big) \Big( \f {k}{p} \Big) =(-1)^{\f {p-1}2} \Big( \f
{k}{p} \Big) = -\Big( \f {k}{p} \Big)$$ where in the last equality we have
used that $p\equiv 3 (4)$.

Taking into account this observation, we obtain the expressions
(\ref{etainvariants}) in the proposition.
\end{proof}

We now look at the simplest case when $p=3$, first considered in
\cite{Pf}. We have $r=0$ and  $\vep_1=(1,1,1,-x(\tf {\pi} 3))=(1,1,1,x(\tf
{\pi} 3 + \pi)),\, \vep_2=(1,1,-1,x(\tf {\pi} 3))$ (in the notation of
Section 4, see (\ref{spinstructure})).
Since $\big( \tf {1}{p} \big) = 1$ we obtain
\begin{eqnarray*}
\eta_{\vep_1}(0)&=& \tf {-2}{\sqrt 3} ( \big( \tf {1}{3} \big) \cot (\tf
{\pi}3))= -\tf 2 {\sqrt 3} \cdot \tf 1{\sqrt 3}= -\tf 23,\\
\eta_{\vep_2}(0)&=& \tf {-2}{\sqrt 3} ( -\big( \tf {1}{3} \big)
\text{cosec} (\tf {\pi}3)) = \tf 2 {\sqrt 3} \cdot \tf 2{\sqrt 3} = \tf
{4}{3}.
\end{eqnarray*}
It is reassuring to see that  the values are in coincidence with those in \cite{Pf}, after all these calculations.

\medskip
To conclude this section,  we shall give explicitly the eta
invariants for all $p$-manifolds in the family, $p=4r+3$ prime,
$7\le p\le 503$, obtained with the help of a computer, using
formulas (\ref{etainvariants}), in the case $(\rho,V)=(1, \C)$. We
also give some values of $d_0(\vep_1)$, the dimension  of the
space of harmonic spinors using (\ref{halfharmspinors}) and
(\ref{cosines}). Note that $d_0(\vep_2)=0$ by Theorem \ref{main}.
We summarize the information in the following table:
\medskip

\renewcommand{\arraystretch}{1}
\begin{center}
\begin{tabular}{|c|c|c|c|r|}
\hline  
$r$  & $p=4r+3$ & $\eta_{\vep_1}$  & $\eta_{\vep_2}$ & $d_0(\vep_1)$ \\
\hline  
0  & 3  & -$\tf23$ &$\tf 43$& 0\\
1  & 7  & -2 & 0 & 2 \\
2  & 11 & -2 & 4 & 2 \\
4  & 19 & -2 & 4 & 26 \\
5  & 23 & -6 & 0 & 90 \\
7  & 31 & -6 & 0 & 1058 \\
10 & 43 & -2 & 4 & 48770 \\
11 & 47 & -10& 0 & 178482\\
14 & 59 & -6 & 12 & 9099506 \\
16 & 67 & -2 & 4 & 128207978\\
17 & 71 & -14& 0 & 483939978 \\
\hline
\end{tabular}
\renewcommand{\arraystretch}{1.08}
\hspace{0.5cm}
\begin{tabular}{|c|c|c|c|}
\hline  
$r$  & $p=4r+3$ & $\eta_{\vep_1}$  & $\eta_{\vep_2}$ \\
\hline  
19 & 79 & -10& 0  \\
20 & 83 & -6 & 12 \\
25 & 103& -10& 0  \\
26 & 107& -6 & 12 \\
31 & 127& -10& 0  \\
32 & 131& -10& 20 \\
34 & 139& -6 & 12 \\
37 &151 & -14& 0  \\
40 & 163& -2 & 4  \\
41 & 167& -22 & 0 \\
\hline
\end{tabular}
\end{center}

\vspace{.5cm}
\begin{center}
\begin{tabular}{|c|c|c|c|}
\hline  
$r$  & $p$ & $\eta_{\vep_1}$  & $\eta_{\vep_2}$ \\
\hline  
44 &179 & -10& 20 \\
47 &191 & -26& 0 \\
49 & 199& -18& 0  \\
52 & 211& -6 & 12 \\
55 & 223& -14& 0  \\
56 & 227& -10& 20 \\
59 & 239& -6 & 12 \\
62 &251 & -30& 28  \\
65 & 263& -26& 0  \\
67 & 271& -22 & 0 \\
\hline
\end{tabular}
\hspace{.475cm}
\begin{tabular}{|c|c|c|c|}
\hline  
$r$  & $p$ & $\eta_{\vep_1}$  & $\eta_{\vep_2}$ \\
\hline  
70 & 283& -6 & 12   \\
76 & 307& -6 & 12   \\
77 & 311& -38& 0    \\
82 & 331& -6 & 12   \\
86 & 347& -10& 20   \\
89 & 359& -38& 0    \\
91 & 367& -18& 0    \\
94 &379 & -6 & 12   \\
95 & 383& -34& 0    \\
104& 419& -18 & 36  \\
\hline
\end{tabular}
\hspace{.475cm}
\begin{tabular}{|c|c|c|c|}
\hline  
$r$  & $p$ & $\eta_{\vep_1}$  & $\eta_{\vep_2}$ \\
\hline  
107& 431& -42& 0  \\
109& 439& -30& 0  \\
110& 443& -10& 20 \\
115& 463& -14& 0  \\
116& 467& -14& 28 \\
119& 479& -50& 0  \\
121& 487& -14& 0  \\
122&491 & -18& 36 \\
124& 499& -6 & 12 \\
125& 503& -42 & 0 \\
\hline
\end{tabular}
\end{center}


\section{APPENDIX: some facts on spin groups and spin representations}

In this appendix we collect some facts on conjugacy classes on
$\text{Spin}(n)$ and on spin representations that are used in the
body of the paper. For standard facts on spin geometry we refer to
\cite{LM} or \cite{Fr2}.

\msk
We consider  $(L, \text{S})$, an irreducible complex
representation  of the Clifford algebra $\C l(n)$, restricted to
$\text{Spin}(n)$. The complex vector space $\text{S}$ has
dimension $2^m$ with $m=[\tfrac n2]$ and is called the spinor
space. We have that $\text{S} = \sum_{I\subset\{1,\ldots,m\}}
\text{S}_{\lambda_I}$ where $\text{S}_{\lambda_I}$ denotes the
weight space corresponding to the weight
\begin{equation} \label{e.weights}
\lambda_I=\f{1}{2}\Big(\sum_{i=1}^m \varepsilon_i\Big) -\sum_{i\in
I}\varepsilon_i.
\end{equation}
Here  
$\vep_j$ is given on the Lie algebra of $T$, $\mathfrak{t}$, by
$\vep_j( \sum_{k=1}^m c_k e_{2k-1}e_{2k})= 2ic_j$. All weights
have multiplicity 1. If $n$ is odd, then $(L,\s)$ is irreducible
for $\text{Spin}(n)$ and is called the {\em spin representation}.
If $n$ is even, then the subspaces
\begin{equation}
\text{S}^+ := \sum_{|I|=even} \text{S}_{\lambda_I},\qquad \qquad
\text{S}^- := \sum_{|I|=odd} \text{S}_{\lambda_I}.
\end{equation}
are $\text{Spin}(n)$-invariant and irreducible of dimension $2^{m-1}$. If
$L^\pm$ denote the restricted action of $L$ on $\s^\pm$ then
$(L^\pm,\s^{\pm})$ are called the {\em half-spin representations} of
$\text{Spin}(n)$. We shall write $(L_n,\text{S}_n)$ (resp.\@
$(L_n^\pm,\s_n^{\pm})$) for $(L,\text{S})$, (resp.\@ $(L^\pm,\s^{\pm})$),
when we wish to specify the dimension.
We have the following wellknown facts:
\begin{equation}\label{spinrestrictions}
\begin{array}{ll} ({L_n ^\pm}_{|\text{Spin}(n-1)}, \s_n^\pm) \simeq (L_{n-1},\s_{n-1})
& \quad \text{if $n$ is even,} \\
({L_n}_{|\text{Spin}(n-1)}, \s_n) \simeq (L_{n-1}^+,\s_{n-1}^+)\oplus
(L_{n-1}^-,\s_{n-1}^-) & \quad \text{if $n$ is odd.}
\end{array}
\end{equation}

The next lemma gives the values of the characters of $L_n$ and ${L_n^\pm}$ on
elements of $T$.
\begin{lema}\label{spincharacters}
If $n=2m$, then
\begin{equation}\label{spinpmchars}
\chi_{_{L^\pm_n}} (x(t_1,\dots,t_m))= 2^{m-1}\Big( \prod_{j=1}^m \cos t_j
\pm i^m \prod_{j=1}^m  \sin t_j \Big).
\end{equation}
If $n=2m$ or $n=2m+1$, then
\begin{equation}\label{cosines}
\chi_{_{L_n}} (x(t_1,\dots,t_m))= 2^{m} \prod_{j=1}^m \cos t_j .
\end{equation}
\end{lema}
\begin{proof}
Assume first that $n=2m$ is even and proceed by induction on $m$.
For $m=1$, (\ref{spinpmchars}) clearly holds. Assume it holds for
$n=2m$. Set $I_m = \{1,\ldots,m\}$. Now, $\chi_{_{L^+_{n+2}}}
(x(t_1,\dots,t_{m+1}))$ equals
\begin{eqnarray*}
&=& e^{i(\sum_1^{m+1} t_j)}
\sum_{{\scriptsize \begin{array}{c} I\subset I_{m+1} \\ |I| \mbox{ even} \end{array} }} e^{-2i \sum_{j \in I} t_j} \\
&=& e^{i t_{m+1}}  e^{i(\sum_1^m t_j)} \; \Big(\!\!\sum_{{\scriptsize
  \begin{array}{c} I\subset I_m\\|I| \mbox{ even} \end{array}  } } e^{-2i \sum_{j \in I} t_j}  \;+\; e^{-2i t_{m+1}}\!\!\sum_{{\scriptsize
   \begin{array}{c} I\subset I_m\\|I| \mbox{ odd} \end{array} } } e^{-2i \sum_{j \in I} t_j}\Big) \\
&=& e^{i t_{m+1}} \,\chi_{_{L^+_n}} \big(x(t_1,\dots,t_{m})\big)
\,
+\, e^{-i t_{m+1}} \,\chi_{_{L^-_n}}\big(x(t_1,\dots,t_m)\big) \\
&=& 2^{m-1} \, \Big( e^{i t_{m+1}} \big(\prod_{j=1}^m \cos t_j +
i^m \prod_{j=1}^m  \sin t_j\big) +  e^{-i t_{m+1}}
\big(\prod_{j=1}^m \cos t_j - i^m
\prod_{j=1}^m  \sin t_j\big) \Big)\\
&=& 2^{m} \Big(\prod_{j=1}^{m+1} \cos t_j + i^{m+1}
\prod_{j=1}^{m+1} \sin t_j \Big).
\end{eqnarray*}
The calculation for $\chi_{_{L^-_{n+2}}}$ is analogous. By adding
$\chi_{_{L^+_n}}(x(t_1,\dots,t_{m}))$ and
$\chi_{_{L^-_n}}(x(t_1,\dots,t_{m}))$ we get the asserted expression for
$\chi_{_{L_n}}(x(t_1,\dots,t_{m}))$ if $n=2m$. If $n=2m+1$, then
$\chi_{_{L_n}}(x(t_1,\dots,t_{m})) =
\chi_{_{L_{n-1}}}(x(t_1,\dots,t_{m}))$, hence the result follows.
\end{proof}

The next lemma gives some useful facts on conjugacy classes of elements in $\text{Spin}(n)$. We include a proof for completeness.

\begin{lema}\label{conjugacy}
Let $x,y \in \text{Spin}(n-1)$ be conjugate in $\text{Spin}(n)$.

\noindent (i) If $n$ is even, then $x,y$ are conjugate in
$\text{Spin}(n-1)$.

\noindent (ii) If $n$ is odd, then $y$ is conjugate to $x$ or to $-e_1 x
e_1$ in $\text{Spin}(n-1)$.
\end{lema}
\begin{proof}
If $n=2m$ is even, the restriction map from the representation ring
$R(\text{Spin}(2m))$ to $R(\text{Spin}(2m-1))$ is onto, hence  the
assertion in the lemma follows.

If $n=2m+1$, we may assume that $x=x(t_1, \dots, t_m), y=x(t'_1, \dots,
t'_m)$ lie in the maximal torus $T$, where $x(t_1, \dots, t_m)=
\prod_{j=1}^m (\cos t_j + \sin t_j e_{2j-1}e_{2j})$.

Now, if $x$ and $y$  are conjugate in $\text{Spin}(2m+1)$, then
$\mu(x), \mu(y)$ are conjugate in $\text{SO}(2m+1)$ and this
implies that,  after reordering, we must have $t'_i=\pm t_i$, for
$1\le i \le m$.

Furthermore if $1\le j\le m$ we have
\begin{eqnarray*}
e_{2j-1}e_n x(t_1,\ldots,t_m) (e_{2j-1}e_n)^{-1} &=&
e_{2j-1}x(t_1,\ldots,t_j,\ldots,t_m) (e_{2j-1})^{-1}\\  &=&
x(t_1,\ldots,-t_j,\ldots,t_m).
\end{eqnarray*}
 Hence, if $x= x(t_1,\ldots,t_m)$, then
\begin{equation*}
(e_{2j-1} e_{2k-1})\, x \, (e_{2j-1} e_{2k-1})^{-1}=
x(t_1,\ldots,-t_j,\ldots,-t_k, \ldots,t_m)
\end{equation*}
for $1 \le j,k \le m$. Thus, for fixed $t_1,\ldots,t_m$, among the
elements of the form $x(\pm t_1, \dots, \pm t_m)$, there are at
most two conjugacy classes in $\text{Spin}(2m)$  represented by
$x(\pm t_1, t_2 \dots,  t_m)$ and $x(- t_1, t_2 \dots, t_m)= -e_1
x(t_1, t_2 \dots,  t_m) e_1$.

Now by Lemma \ref{spincharacters}, we have that
\begin{equation*}
\chi_{_{L^\pm_{n-1}}} (x(t_1,\dots,t_m))= 2^{m-1} \Big(
\prod_{j=1}^m \cos t_j \pm i^m \prod_{j=1}^m  \sin t_j \Big).
\end{equation*}
This implies that $x(t_1, t_2 \ldots,  t_m)$ and $x(- t_1, t_2 \ldots,
t_m)$ are not conjugate unless $t_j \in \pi \Z$ for some $j$. On the other
hand, if this is the case, then clearly $e_1e_{2j-1}\in \text{Spin}(n-1)$
conjugates one element into the other. This completes the proof of the
lemma.
\end{proof}

\begin{rem} The lemma shows that generically, if $n$ is odd, $x(t_1, t_2 \dots,  t_m)$ and
$x(-t_1, t_2 \dots,  t_m)$ are conjugate in $\text{Spin}(n)$ but  not  in $\text{Spin}(n-1)$.
\end{rem}

We now consider the special case when  $t_i \in \frac {\pi} 2 \Z$ for all
$i$,
then $\mu(x)$ has order 2 (or 1). Set $g_h = e_1e_2\dots e_{2h-1}e_{2h}
\in \spin$ for $1\le h \le m= [\frac n2]$. Thus $g_h =
x(\underbrace{\tfrac \pi 2,\dots,\tfrac \pi 2}_h,0, \dots, 0)$ and $-g_h =
x(\underbrace{-\tfrac \pi 2,\tfrac \pi 2, \dots,\tfrac \pi 2}_h,0, \dots,
0)$ .

\begin{coro}\label{tracespin}If $h<m$, then $g_h$ and $-g_h$ are conjugate in
$\text{Spin}(n-1)$. If $h=m$ and $n=2m$, then
$\chi_{_{L^\pm_n}}(g_m)= \pm 2^{m-1}i^m$, hence $g_m$ and $-g_m$
are not conjugate. If $h=m$ and $n=2m+1$,  then
$\chi_{L_{n-1}^+}(\pm g_m)=\pm 2^{m-1}i^m $, hence $g_m$ and
$-g_m$ are conjugate in $\text{Spin}(n)$ but not in
$\text{Spin}(n-1)$.
\end{coro}
\begin{proof}
The first assertion in the corollary follows immediately from the
proof of Lemma \ref{conjugacy}. The remaining assertions are clear
in light of Lemma \ref{spincharacters}. Indeed, for $h<m$ we have
$e_1e_n g_h (e_1 e_n)^{-1} =-g_h$.
\end{proof}

Recall that for any $u\in \R^n \smallsetminus \{0 \}$, left
Clifford multiplication by $u$ on $\s$ is given by $u \cdot w =
L(u)(w)$ for $w\in \text{S}$. We fix $\langle\, ,\, \rangle$ an
inner product on $\s$ such that $L(u)$ is skew Hermitian, hence
$\langle\, ,\, \rangle$ is $Spin(n)$-invariant. Note that $L(u)^2
=-\|u\|^2\I$. Hence, $\s$ decomposes $\text{S}=\text{S}_u^+\oplus
\text{S}_{u}^-$, where $\text{S}^\pm_u$ denote the eigenspaces, of
dimension $2^{m-1}$, of $L(u)$ with eigenvalues $\mp i\|u\|$.

\begin{defi} If $u \in \R^n \smallsetminus \{0 \}$ set
\begin{equation}\label{spinu}
\text{Spin}(n-1,u):= \{g \in \text{Spin}(n) : gug^{-1}=u \}.
\end{equation}
\end{defi}
Clearly $\text{Spin}(n-1,e_n)=\text{Spin}(n-1)$ and for general $u$, if
$h_u \in \text{Spin}(n)$ is such that $h_u u {h_u}^{-1} = \|u\| e_n$, then
$h_u \text{Spin}(n-1,u) h_u^{-1} = \text{Spin}(n-1)$. We note that for any
$g \in \text{Spin}(n-1,u)$, $L(g)$ commutes with $L(u)$, hence  $L(g)$
preserves the eigenspaces $\text{S}^\pm_u$.

The following lemma is used in the proof of Theorem \ref{main}. 
\begin{lema}\label{lemaSv}
Let $\text{Spin}(n-1,u)$ be as in (\ref{spinu}). Then as
$\text{Spin}(n-1)$-modules: $\s_{e_n}^\pm \simeq (L_{n-1}^\pm,
\s_{n-1}^\pm)$ if $n$ is odd and $\s_{e_n}^\pm \simeq (L_{n-1},\s_{n-1})$
if $n$ is even. As $\text{Spin}(n-1,u)$-modules we have that $\s_u^\pm =
L(h_u) \s_{n-1}^\pm$, if $n$ is odd, and $\s_u^\pm = L(h_u) \s_{n-1}$, if
$n$ is even,  with action given by $L(h_u x h_u^{-1}) =
L(h_u)L(x)L(h_u^{-1})$ for any $x \in \text{Spin}(n-1)$.
\end{lema}
\begin{proof}
$L(e_n)$ commutes with the action of $\text{Spin}(n-1)$ on $\s_n$ and, on
the other hand, $\s_n= \s_{n-1}^+ \oplus \s_{n-1}^-$ as a
$\text{Spin}(n-1)$-module.

If $n$ is odd, then $\s_{n-1}^\pm$ are inequivalent representations of
Spin$(n-1)$, hence $L(e_n)\s_{n-1}^\pm = \s_{n-1}^\pm$
and by Schur's lemma, $L(e_n)$ must act by multiplication by a scalar on  each of them. 
By using the explicit description of $L$ in \cite{Kn}, p.\@ 286--288, one
verifies that $L(e_n)$ acts by $\mp i$ on $\s_{n-1}^\pm$, that is
$\s_{e_n}^\pm \simeq \s_{n-1}^\pm$.

If $n$ is even, then
$\s_{n}^\pm$ both restrict to $\s_{n-1}$ as $\text{Spin}(n-1)$-modules. Since the $\pm i$-eigenspaces of $L(e_n)$ are stable by  $\text{Spin}(n-1)$, they must both be equivalent to  $\s_{n-1}$. 

The remaining assertions  are easily verified.
\end{proof}


\end{document}